\documentclass[12pt]{article}
\usepackage{graphicx}
\usepackage{comment}
\usepackage{amssymb,amsmath,amsthm}
\usepackage[english]{babel}
\usepackage[usenames]{color}
\usepackage{xcolor}
\usepackage{authblk}
\usepackage{tikz}
\usetikzlibrary{arrows.meta}
\usepackage{xxcolor}
\usepackage{array}
\usepackage{type1cm}
\usepackage{soul}
\def\N{{\rm I\kern-.20em N}}
\def\R{{\rm I\kern-.17em R}}
\def\C{{\rm I\kern-.17em C}}

\def\bkC{{\rm \kern.24em\vrule
width.05em height1.4ex depth-.05ex\kern-.26emC}}
\definecolor{laranja}{rgb}{1,0.5,0}

\newcommand{\g}[1]{g_{#1}^{\alpha}}
\newcommand{\z}[1]{z_{#1}^{\alpha}}

\begin{document}
\title{Fractional differential problems with numerical anti-reflective boundary conditions: a computational/precision analysis and numerical results}

\author[1]{Ercília Sousa}
\author[2]{Cristina~Tablino-Possio}
\author[3]{Rolf Krause}
\author[4,5]{Stefano~Serra-Capizzano}

\affil[1]{{\small Department of Mathematics, University of Coimbra, Coimbra, Portugal (email: ecs@mat.uc.pt).}}
\affil[2]{{\small Department of Mathematics and Applications, University of Milano - Bicocca, Milano, Italy (email: cristina.tablinopossio@unimib.it).}}
\affil[3]{{\small Center for Computational Medicine in Cardiology, University of Italian Switzerland, Lugano, Switzerland (email: rolf.krause@usi.ch).}}
\affil[4]{{\small Department of Science and High Technology, University of Insubria, Como, Italy (email: s.serracapizzano@uninsubria.it).}}
\affil[5]{{\small Department of Information Technology, Uppsala University, Uppsala, Sweden (email: serra@it.uu.se).}}
\maketitle
%-------------------------------------------------------------------------------------------------------
\begin{abstract}
 Twenty years ago the anti-reflective numerical boundary conditions (BCs) were introduced in a context of signal processing and imaging, for increasing the quality of the reconstruction of a blurred signal/image contaminated by noise and for reducing the overall complexity to that of few fast sine transforms i.e. to $O(N\log N)$ real arithmetic operations, where $N$ is the number of pixels.
Here for quality of reconstruction we mean a better accuracy and the elimination of boundary artifacts, called ringing effects.
Now we propose numerical anti-reflective BCs in the context of nonlocal problems of fractional differential type: the goals are the same i.e. a smaller approximation error and the reduction of boundary artifacts.
 In the latter setting, we compare various types of numerical BCs, including the anti-symmetric ones considered in the case of fractional differential problems for modeling reasons. More in detail, given important similarities between anti-symmetric and anti-reflective BCs, we compare them from the perspective of computational efficiency, by considering nontruncated and truncated versions and also other standard numerical BCs such as reflective/Neumann. Several numerical tests, tables, and visualizations are provided and critically discussed. The conclusion is that the truncated numerical anti-reflective BCs perform better, both in terms of accuracy and low computational cost.
\end{abstract}

%\tableofcontents
%-------------------------------------------------------------------------------------------------------
\section{Introduction}\label{sec:intro}
Twenty years ago the anti-reflective boundary conditions (BCs) were introduced \cite{S-AR-proposal} in a context of signal processing and imaging, for increasing the quality of the reconstruction of a blurred signal/image contaminated by noise and for reducing the overall complexity to that of few fast sine transforms i.e. to $O(N\log N)$ real arithmetic operations, where $N$ is the number of pixels: see also \cite{paperAR,paperART,imaging-AR-2,imaging-AR-1,imaging-AR-0} for a series of results regarding the anti-reflective algebra of matrices, the low complexity transform, and deblurring techniques associated with anti-reflective BCs. In other words, the anti-reflective BCs represent a numerical trick for obtaining higher precision at low computational cost. To be more specific for quality of reconstruction we mean a better accuracy and the elimination of disturbing boundary artifacts, called ringing effects \cite{artifact 1,artifact 2,artifact 3}, also in connection with fractional operators \cite{artifact - fract}.
\par
Here we propose the anti-reflective BCs in the context of nonlocal problems of fractional differential type.
The motivation is twofold:
\begin{itemize}
\item the link between fractional models and imaging problems relies on the nonlocality of the underlying operator;
\item rough artifacts stemming from the boundaries have been observed also in a fractional differential setting \cite{BU2014,rough-artifact}.
\end{itemize}
Indeed, regarding the first item, in the imaging case an integral operator describes the model, while the fractional operators are also described in terms of a global integral and here they share their nonlocal nature.
The difference is that fractional equations are endowed with physical BCs, while this is not the case in general in signal and image processing. Hence in this work we have to be careful in making a distinction between physical and numerical BCs. Interestingly enough, since the observation in \cite{BU2014}, it is now generally accepted that also in the numerical approximation of fractional differential equations boundary artifacts appear and they become worse and worse, when increasing the matrix-size $N$ (see also \cite{rough-artifact}).  Hence the goal of setting additional numerical BCs for diminishing the boundary artifacts and for having an optimal $O(N\log N)$ complexity is of concrete interest in real-world applications modeled by fractional operators (see e.g. \cite{fract-appl1,fract-appl2} and references therein). We also remark that fractional operators are now used also in imaging and hence making a connection between the two areas is crucial for a mutual benefit.

We consider a specific class of fractional equations, but the idea is general and can be applied in a much wider context.
We start to define  the left Riemann-Liouville fractional derivative of order
$\alpha$, with $1<\alpha<2$ and $x\in \R$, as
$$
{}_{-\infty}^{RL}D_{x}^\alpha u(x,t) = \frac{1}{\Gamma(2-\alpha)}\frac{\partial^2}{\partial x^2}
\int_{-\infty}^{x} {u(\xi,t)}{(x-\xi)^{1-\alpha}}d\xi,
$$
and the right Riemann-Liouville fractional derivative of order $\alpha$, with $1<\alpha<2$ and $x\in \R$, expressed as
$$
{}_{x}^{RL}D_{\infty}^\alpha u(x,t) = \frac{1}{\Gamma(2-\alpha)}\frac{\partial^2}{\partial x^2}
\int_{x}^{\infty} {u(\xi,t)}{(\xi-x)^{1-\alpha}}d\xi.
$$
Setting $ \kappa_\alpha = ( \cos(\pi\alpha/2) )^{-1} $, the following class of general weighted operators is considered
$$
\Delta_\beta^{\alpha/2} u(x,t) = \frac{1+\beta}{2}\ {}_{-\infty}^{RL}D_{x}^\alpha u(x,t)
+  \frac{1-\beta}{2} \ {}_{x}^{RL}D_{\infty}^\alpha u(x,t),
$$
with $\alpha\in (1,2)$, $\beta\in (-1,1)$ so that $\Delta_\beta^{\alpha/2}$ becomes a linear convex combination of the given left and right derivatives, with $\Delta_0^{\alpha/2}\equiv \Delta^{\alpha/2}$ being their arithmetic mean for $\beta=0$ and
\begin{equation}\label{beta=0}
(-\Delta)^{\alpha/2} u(x,t) = \frac{1}{2\cos(\pi\alpha/2)} \left[
{}_{-\infty}^{RL}D_{x}^\alpha u(x,t)+{}_{x}^{RL}D_{\infty}^\alpha u(x,t)
\right],
\end{equation}
according to \cite[p.11]{lis2020}.

Regarding further connections,  we recall the following class of nonlinear fractional problems appearing in several works \cite{dip2023,dip2023p2,zhu2022},
usually in $d$-dimensional domains and without the time variable. The precise formulation can be described by the following two equations
\begin{align}
&(-\Delta)^{\alpha/2} u(x) = f(u), \quad x \in \R^d,
\label{riesz}\\
& u(x',-x_d)  =  - u(x',x_d), \quad  x=(x',x_d)  \in \R^d,
\label{reflective1_initial}
\end{align}
with $x'=(x_1,\dots,x_{d-1})$, $(-\Delta)^{\alpha/2}$ being the fractional Laplacian with fractional order $\alpha\in (1,2)$, and
$f(u)=-c(x)u^p$, $c,p\geq 0$.

When considering $d=1$, condition (\ref{reflective1_initial}) reduces to $u(-x)=-u(x)$ with $x_d\equiv x$ and $x'$ not present, while the operator in (\ref{riesz}) in $\R$ is defined as in \cite[p.11]{lis2020} with the precise expression reported in (\ref{beta=0}), if we omit the time variable. We observe that the present physical anti-symmetry is somehow a special case of the numerical anti-reflection used in imaging: the latter statement is crystal clear by comparing equations (\ref{structure-anti-symm}), (\ref{structure-anti-refl}) with \cite{S-AR-proposal}[Eq. (3.3)], since the matrix structures emerging in the present work belong to the anti-reflective algebra introduced in \cite{S-AR-proposal}. We stress again and it is worth noticing that the conditions (\ref{reflective1_initial}) are called anti-symmetric BCs and characterize the continuous equation for modeling reasons, while numerical anti-reflective BCs are not part of a continuous model as the anti-symmetric BCs in (\ref{reflective1_initial}).

%as (\ref{riesz})

Therefore, inspired by the previous considerations and by the connections emphasized in the previous lines, in the current work we consider
the one-dimensional  time-dependent problem, given by
\begin{align*}
\frac{\partial u}{\partial t} (x,t) &=  \kappa_\alpha \Delta_\beta^{\alpha/2} u(x,t),\quad x \in (a,b), \ t>0,\\
\end{align*}
$\alpha\in (1,2)$, $\beta\in (-1,1)$, with initial condition in time and when we reconstruct the behavior of the solution outside the interval $[a,b]$ by making use of numerical zero Dirichlet BCs, numerical reflection, anti-symmetry and anti-reflection.

The present work is organized as follows. In Section \ref{sec:formulation} we describe the problem in an open domain, while in Sections \ref{sec:anti-sym} and \ref{sec:anti-refl} we enforce the two types of numerical BCs connected with the presence of walls and we discuss the implication of the infinite number of (Fourier) coefficients coming from the nonlocal nature of the underlying operators. Section \ref{sec:num} deals with several experiments and visualizations, which look very interesting from a numerical viewpoint due to a very low complexity of $O(N\log N)$ real arithmetic operations, with $N$ being the number of space grid points: as expected from the theoretical study in \cite{S-AR-proposal} in the context of signal processing and imaging, the numerical precision is higher when considering the anti-reflective choice with respect to all the other considered numerical BCs, even in our considered fractional setting. Section \ref{sec:truncations} is devoted to the study of the proper truncations for imposing that the resulting matrices lie in the anti-reflective matrix algebra, so further diminishing the computational cost of $O(N\log N)$ real arithmetic operations, where the hidden constant is substantially lower with respect to the previous nontruncated case. Several numerical results are reported and discussed for checking the computational and precision advantages of our proposals. Finally, Section \ref{sec:end} contains final remarks and a mention to a list of open problems.
%-------------------------------------------------------------------------------------------------------
\section{Open domain}\label{sec:formulation}
 Suppose first that we have the {time dependent} fractional diffusion equation in the open space domain i.e.
\begin{equation*}
 \frac{\partial u}{\partial t}(x,t)
   =  \kappa_\alpha  \Delta_\beta^{\alpha/2} u(x,t), \quad x \in \R, \ t>0,
 \end{equation*}
{with
$$
\Delta_\beta^{\alpha/2} u(x,t) =  \frac{1+\beta}{2}\ {}_{-\infty}^{RL}D_{x}^\alpha u(x,t) + \frac{1-\beta}{2} \ {}_{x}^{RL}D_{\infty}^\alpha u(x,t),
$$
with $\alpha\in (1,2)$, $\beta\in (-1,1)$ so that $\Delta_\beta^{\alpha/2}$ becomes a linear convex combination of the given left and right derivatives, with $\Delta_0^{\alpha/2}\equiv \Delta^{\alpha/2}$ being their arithmetic mean for $\beta=0$.
\par
Hereafter, we will consider the Gr{\"u}nwald-Letnikov approximations  of the left and right Riemann-Liouville  fractional derivatives  at $(x_j,t_n)$, that is  \begin{align}
{}_{-\infty}^{RL}D_{x}^\alpha u (x_j,t_n) & \approx \frac{1}{(\Delta x)^\alpha}\sum_{k=0}^{\infty} g_{k}^\alpha u(x_{j+1-k},t_n), \label{glapproxl} \\
{}_{x}^{RL}D_{\infty}^\alpha u (x_j,t_n) & \approx \frac{1}{(\Delta x)^\alpha}\sum_{k=0}^{\infty} g_{k}^\alpha u(x_{j-1+k},t_n), \label{glapproxr}
\end{align}
where the Gr{\"u}nwald-Letnikov coefficients are defined trought the subsequent recurrence formula for all $\alpha>0$ as
\begin{equation*}
g_0^\alpha=1, \quad g_{k+1}^\alpha = -\frac{\alpha-k}{k+1} g_{k}^\alpha, \quad k\geq 0.
%\label{gk}
\end{equation*}
Let $U_j^n$ represent the approximate solution of $u(x_j,t_n)$ in the discrete domain and let ${\bf U}^{n}=[ U_{-N}^n, \dots, U_{N}^n]^T$ taking
into consideration  that the function goes to zero as we go to infinity.
\par
Then the $\theta-$method for the time integration in matrix form gives
\[
\frac{{\bf U}^{n+1} - {\bf U}^{n}}{\Delta t} = \frac{\kappa_\alpha}{(\Delta x)^\alpha} A_\beta \left ( \theta {\bf U}^{n+1} + (1-\theta) {\bf U}^{n} \right ),
\]
or, equivalently,
\[
(I -  \mu_\alpha \theta A_\beta) {\bf U}^{n+1} = (I +  \mu_\alpha (1-\theta) A_\beta)  {\bf U}^{n},
\]
where the matrix $A_\beta/{(\Delta x)^\alpha}$ represents the chosen approximation of $\Delta_\beta^{\alpha/2} (\cdot)$, $I$ represents the identity matrix, and
$\mu_\alpha={\kappa_\alpha \Delta t}/{(\Delta x)^\alpha}$.  \\
As well known, the Esplicit Euler method is obtained for $\theta =0$, the Implicit Euler method for $\theta =1$, and the Crank-Nicolson method for $\theta=1/2$.
Clearly, the application of the Explicit Euler method simply requires matrix-vector multiplications with Toeplitz structures.
In Section \ref{sec:num} we will consider the more interesting case of application of both Implicit Euler  and Crank-Nicolson methods where we have to solve the structured linear systems above. These methods with $\theta =1, 1/2$ have much better stability and approximation properties.
\par
Thus, to form the matrix $A_\beta$, we simply need to
% where $\mu_\alpha={\kappa_\alpha \Delta t}/{(\Delta x)^\alpha}$.
consider \eqref{glapproxl}-\eqref{glapproxr} and the assumption ${\bf U}^{n}=[ U_{-N}^n, \dots, U_{N}^n]^T$, so that }
$$
A_\beta = \frac{1+\beta}{2} A_L+\frac{1-\beta}{2} A_R,
$$
$$
A_L = \left[
\begin{array}{cccccc}
g_1^\alpha  & g_0^\alpha & 0 & \dots & 0 &0\\
g_2^\alpha & g_1^\alpha & g_0^\alpha& \dots & 0 &0 \\
g_3^\alpha  & g_2^\alpha & g_1^\alpha & \dots & 0 & 0 \\
\vdots & \vdots& \vdots &&\vdots&\vdots\\
g_{2N+1}^\alpha & g_{2N}^\alpha & g_{2N-1}^\alpha & \dots & g_2^\alpha & g_1^\alpha
\end{array}
\right],
$$
$A_R = A_L^T,$ that is
$$
A_R = \left[
\begin{array}{cccccc}
g_1^\alpha  & g_2^\alpha & g_3^\alpha & \dots & g_{2N}^\alpha &g_{2N+1}^\alpha\\
g_0^\alpha & g_1^\alpha & g_2^\alpha& \dots & g_{2N-1}^\alpha & g_{2N}^\alpha \\
0  & g_0^\alpha & g_1^\alpha & \dots &  & g_{2N-1}^\alpha \\
\vdots & \vdots& \vdots &&\vdots&\vdots\\
0 & 0 & 0 & \dots & g_0^\alpha & g_1^\alpha
\end{array}
\right].
$$
{It is evident} both $A_L$ and $A_R$ and hence $A_\beta$ share a Toeplitz structure. Furthermore, the choice of $\beta=0$ leads to the Riesz operator (fractional Laplacian in 1D) and
$$
A_0 = \frac{1}{2} A_L+\frac{1}{2} A_R.
$$
Therefore, in accordance with the continuous operator, we find a symmetric Toeplitz matrix of the form
$$
A_0 = \left[
\begin{array}{cccccc}
2g_1^\alpha  & g_0^\alpha +g_2^\alpha & g_3^\alpha  & \dots &  g_{2N}^\alpha &g_{2N+1}^\alpha \\
g_2^\alpha+g_0^\alpha & 2g_1^\alpha & g_0^\alpha+g_2^\alpha & \dots & g_{2N-1}^\alpha & g_{2N}^\alpha \\
g_3^\alpha  & g_2^\alpha +g_0^\alpha& 2g_1^\alpha & \dots & g_{2N-2}^\alpha& g_{2N-1}^\alpha  \\
\vdots & \vdots& \vdots &&\vdots&\vdots\\
g_{2N+1}^\alpha & g_{2N}^\alpha & g_{2N-1}^\alpha & \dots & g_2^\alpha +g_0^\alpha& 2g_1^\alpha
\end{array}
\right].
$$
\par
{ Notice that we could even  consider} other types of approximations that would end up having different values for the coefficients. However, when using uniform griddings, the resulting matrix structure
would remain unchanged, since it is essentially decided by the nonlocal nature of the underlying continuous operator. Hence, independently of the approximation scheme as long as the grid points are equispaced, the resulting matrix structure is the same and essentially of dense Toeplitz type (i.e. shift invariant in the imaging terminology \cite{hansen-nagy-o}).
We point out that the latter represents the main target of our present study, since the computational cost of the related algorithms strongly depends on matrix structural features.
Regarding structured matrices, in detail a square matrix $A$ of size $m$ is called Toeplitz if ${a}_{j,k}=\alpha_{j-k}$, $\alpha_s \in \mathbb{C}$, $s=1-m,\ldots,m-1$. A square matrix $B$ of size $m$ is called Hankel if ${a}_{j,k}=\beta_{j+k-1}$, $\beta_s \in \mathbb{C}$, $s=1,\ldots,2m-1$. We point out that any Hankel matrix $B$ can be written as $JA$ and $\tilde A J$ where $A,\tilde A$ are Toeplitz and $J$ is the antidiagonal or flip matrix with
$J_{j,k}=1$ if $j+k-1=m$ and zero otherwise. Hence $J$ is a special Hankel matrix.

%-------------------------------------------------------------------------------------------------------
\section{Anti-symmetric numerical BCs}\label{sec:anti-sym}
We move now to the problem with an initial condition $u(x,0)=u_0(x), \ x\in [a,b]$, and numerical BCs outside $(a,b)$
\begin{align}
& \frac{\partial u}{\partial t}(x,t) \label{eq:pb}
   =   \kappa_\alpha \Delta_\beta^{\alpha/2} u(x,t), \quad x\in (a,b).
   \end{align}
More specifically, for the use of values outside the interval $[a,b]$, in analogy with what is done in imaging/signal processing outside the field of values \cite{hansen-nagy-o}, we start by using numerical anti-symmetric BCs
\begin{align}
& u(-x+a,t)  =  -u(x+a,t), \quad \mbox{for all} \quad x < a,
\label{reflective1}
\\
& u(x+b,t)  =  -u(-x+b,t), \quad \mbox{for all} \quad x > b,
\label{reflective2}
\end{align}
 and a equispaced discrete domain $x_j= a+j\Delta x$, $j=0,1, \dots, N$, $x_0=a$, $x_{N}=b$.
As a consequence we find
\begin{align*}
u(-x_j+a,t)  &=  -u(x_j+a,t), \quad \mbox{for all} \quad x < a,\\
u(x_j+b,t)  &=  -u(-x_j+b,t), \quad \mbox{for all} \quad x > b.
\end{align*}
When a numerical anti-symmetric  boundary condition is imposed at $x=a$, from (\ref{reflective1}) we deduce
$$U_{j+1-k}^n=-U_{-j-1+k}^n.
$$
The approximation of the left fractional Riemann-Liouville derivative
becomes
\begin{eqnarray*}
 {}_{-\infty}^{RL}D_{x}^{\alpha,ref} u (x_j,t_n)
 &\approx   &
 {\frac{1}{(\Delta x)^\alpha} \left (\sum_{k=0}^{j+1} g_{k}^\alpha U_{j+1-k}^n
+\sum_{k=j+2}^{\infty} g_{k}^\alpha U_{j+1-k}^n \right ) }\\
 & = &
 {\frac{1}{(\Delta x)^\alpha} \left (\sum_{k=0}^{j+1} g_{k}^\alpha U_{j+1-k}^n
 -\sum_{k=j+2}^{\infty} g_{k}^\alpha U_{k-j-1}^n\right ).}
\end{eqnarray*}
However, the fact that the second sum goes until infinity looks problematic, since the anti-symmetric condition
would send points to the right boundary wall and not inside the interior domain $(a,b)$. As already observed, the interior domain $(a,b)$ is called field of values in imaging and signal processing terminology \cite{hansen-nagy-o}.
In fact, we  decide to stop the second sum as described in what follows, that is,
\begin{align} \label{glapproxl_anti}
 {}_{-\infty}^{RL}D_{x}^{\alpha,ref} u (x_j,t_n)
 &\approx
 %\textcolor{violet}
 { \frac{1}{(\Delta x)^\alpha} \left (\sum_{k=0}^{j+1} g_{k}^\alpha U_{j+1-k}^n
 - \sum_{k=j+2}^{N+j+1} g_{k}^\alpha U_{k-j-1}^n \right ) }
\end{align}
and the latter implies that we consider a bounded wall, for the numerical anti-symmetric boundary,
of the same size as the interior domain.
We observe that the truncation gives an approximation that makes sense if we take into consideration the fact that the
sequence $g_k$ goes to zero as $k$ tends to infinity.

If the  boundary is only at  one side the related jumping problem between boundaries would not exist. With boundaries at both sides, we can suppose that we have this reflecting boundary as an intermediate boundary, where at the final boundary wall the function could be set to zero.

For the numerical anti-symmetric  boundary condition at $x=b$, by following a similar approach,  we infer
\begin{eqnarray*}
 {}_{x}^{RL}D_{\infty}^{\alpha,ref} u (x_j,t_n)
 &\approx &
 {\frac{1}{(\Delta x)^\alpha} \left(\sum_{k=0}^{N-j+1} g_{k}^\alpha U_{j-1+k}^n
 +\sum_{k=N-j+2}^{\infty} g_{k}^\alpha U_{j-1+k}^n\right ).}
 \end{eqnarray*}
We also need to stop at a finite point as previously, that is
\begin{eqnarray*}
 {}_{x}^{RL}D_{\infty}^{\alpha,ref} u (x_j,t_n)
 &\approx &
 { \frac{1}{(\Delta x)^\alpha} \left (\sum_{k=0}^{N-j+1} g_{k}^\alpha U_{j-1+k}^n
 +\sum_{k=N-j+2}^{2N-j+1} g_{k}^\alpha U_{j-1+k}^n \right).}
\end{eqnarray*}
Owing to  $u(x+b) = -u(-x+b)$, for $k\geq N-j+2$, we have
 $$
 U_{j-1+k} = U_{N -N+j-1+k} = - U_{N +N-j+1-k} = -  U_{2N-j+1-k}.
 $$
 Therefore
 \begin{align}\label{glapproxr_anti}
 {}_{x}^{RL}D_{\infty}^{\alpha,ref} u (x_j,t_n)
 & \approx
 { \frac{1}{(\Delta x)^\alpha}\left ( \sum_{k=0}^{N-j+1} g_{k}^\alpha U_{j-1+k}^n
 - \sum_{k=N-j+2}^{2N-j+1} g_{k}^\alpha U_{2N-j+1-k}^n \right ).}
\end{align}
{Thus, by applying again the $\theta-$method, the matrix form of the problem is
\begin{equation*}
%{\bf U}^{n+1} = ( I +\mu_\alpha {\bf U}^{n},
(I -  \mu_\alpha \theta A_{\beta}^\textrm{anti}) {\bf U}^{n+1} = (I +  \mu_\alpha (1-\theta) A_{\beta}^\textrm{anti})  {\bf U}^{n},
\end{equation*}
with ${\bf U}^{n}=[ U_{0}^n, \dots, U_{N}^n]^T$, $I$ being the identity matrix, and $\mu_\alpha={\kappa_\alpha \Delta t}/{(\Delta x)^\alpha}$. By combining \eqref{glapproxl_anti} and \eqref{glapproxr_anti}, the matrix $A_{\beta}^\textrm{anti}$ is
expressed as the $(N+1)\times(N+1)$ matrix }
\begin{equation}\label{eq:Abeta_anti}
A_{\beta}^\textrm{anti} =  \frac{1+\beta}{2} A_L^\textrm{anti} +\frac{1-\beta}{2}  A_R^\textrm{anti}
\end{equation}
with
$$
A_{L}^\textrm{anti} = \left[
\begin{array}{cccccccc}
 g_1^\alpha  & g_0^\alpha   \textcolor{brown}{-g_2^\alpha}&  \textcolor{brown}{-g_3^\alpha}  &  \dots &\textcolor{brown}{-g_{N}^\alpha}
 & \textcolor{brown}{-g_{N+1}^\alpha} \\
 g_2^\alpha & g_1^\alpha  \textcolor{brown}{-g_3^\alpha} & g_0^\alpha \textcolor{brown}{-g_4^\alpha} & \dots & \textcolor{brown}{-g_{N+1}^\alpha} &
\textcolor{brown}{-g_{N+2}^\alpha} \\
g_3^\alpha  & g_2^\alpha  \textcolor{brown}{-g_4^\alpha} & g_1^\alpha \textcolor{brown}{-g_5^\alpha} & \dots &  \textcolor{brown}{-g_{N+2}^\alpha}
&   \textcolor{brown}{-g_{N+3}^\alpha}  \\
\vdots & \vdots   &\vdots &&\vdots&\vdots\\
g_{N+1}^\alpha & g_{N}^\alpha  \textcolor{brown}{-g_{N+2}^\alpha} & g_{N-1}^\alpha \textcolor{brown}{-g_{N+3}^\alpha} & \dots
& g_2^\alpha  \textcolor{brown}{-g_{2N}^\alpha}  \textcolor{red}{-g_{0}^\alpha}&  g_1^\alpha \textcolor{brown}{-g_{2N+1}^\alpha}
\end{array}
\right]
$$
and
$$
A_{R}^\textrm{anti}  = \left[
\begin{array}{cccccc}
g_1^\alpha \textcolor{brown}{-g_{2N+1}^\alpha}  & g_2^\alpha \textcolor{brown}{-g_{2N}^\alpha}
\textcolor{red}{-g_0^\alpha} & g_3^\alpha \textcolor{brown}{-g_{2N-1}^\alpha} & \dots & g_{N}^\alpha \textcolor{brown}{-g_{N+2}^\alpha}&g_{N+1}^\alpha\\
g_0^\alpha  \textcolor{brown}{-g_{2N}^\alpha} & g_1^\alpha  \textcolor{brown}{-g_{2N-1}^\alpha}
& g_2^\alpha  \textcolor{brown}{-g_{2N-2}^\alpha} & \dots & g_{N-1}^\alpha  \textcolor{brown}{-g_{N+1}^\alpha} & g_{N}^\alpha \\
0   \textcolor{brown}{-g_{2N-1}^\alpha} & g_0^\alpha \textcolor{brown}{-g_{2N-2}^\alpha} & g_1^\alpha \textcolor{brown}{-g_{2N-3}^\alpha} & \dots &  & g_{N-1}^\alpha \\
\vdots & \vdots& \vdots &&\vdots&\vdots\\
0 \textcolor{brown}{-g_{N+1}^\alpha}& 0 \textcolor{brown}{-g_{N}^\alpha}& 0 \textcolor{brown}{-g_{N-1}^\alpha}& \dots
& g_0^\alpha \textcolor{brown}{-g_{2}^\alpha}& g_1^\alpha
\end{array}
\right].
$$
We stress that the terms appearing because of the presence of the boundary are highlighted in a different color.
{
We also notice as the presence of the correction term $\g{0}$ in the first and last equation directly originates from the fact that the approximation of
\eqref{glapproxl} and \eqref{glapproxr} starts from $k=0$, that is the approximation across the considered point becomes an approximation across the boundary, so requiring a proper reflection inside.}
{When $\beta=0$ we find
$$
A_{0}^\textrm{anti} =  \frac{1}{2} A_L^\textrm{anti} +\frac{1}{2}A_R^\textrm{anti}
$$
so that $2A_{0}^\textrm{anti}$ equals the matrix
\[ \tiny
\begin{bmatrix}
\begin{array}{c|cccccc|c}
2\g{1} \textcolor{brown}{-\g{2N+1}}  &
\g{0}+\g{2} \textcolor{brown}{-\g{2}} \textcolor{brown}{-\g{2N}} \textcolor{red}{-\g{0}} &
\g{3} \textcolor{brown}{-\g{3}} \textcolor{brown}{-\g{2N-1}} &
\ldots &  \ldots &
\g{N-1} \textcolor{brown}{-\g{N-1}} \textcolor{brown}{-\g{N+3}} &
\g{N} \textcolor{brown}{-\g{N}}\textcolor{brown}{-\g{N+2}}  &
\g{N+1} \textcolor{brown}{-\g{N+1}}\\
\hline
\g{0}+\g{2}\textcolor{brown}{-\g{2N}} &
2\g{1}\textcolor{brown}{-\g{3}} \textcolor{brown}{-\g{2N-1}}     &
\g{0}+\g{2} \textcolor{brown}{-\g{4}} \textcolor{brown}{-\g{2N-2}}&
\g{3}\textcolor{brown}{-\g{5}} \textcolor{brown}{-\g{2N-3}} &  \ldots      & \ldots &
\g{N-1} \textcolor{brown}{-\g{N+1}} \textcolor{brown}{-\g{N+1}}                     &
\g{N} \textcolor{brown}{-\g{N+2}}     \\
\g{3}  \textcolor{brown}{-\g{2N-1}}     &
\g{0}+\g{2}\textcolor{brown}{-\g{4}} \textcolor{brown}{-\g{2N-2}}&
\vdots      & \vdots & \vdots & \vdots &
\g{N-2}\textcolor{brown}{-\g{N+2}}\textcolor{brown}{-\g{N}}                         &
\g{N-1} \textcolor{brown}{-\g{N+3}}    \\
\vdots      & \vdots                               & \vdots      & \vdots & \vdots & \vdots & \vdots                          & \vdots      \\
\g{N-1} \textcolor{brown}{-\g{N+3}}    &
\g{N-2}\textcolor{brown}{-\g{N}} \textcolor{brown}{-\g{N+2}}    &
\vdots                          & \vdots & \vdots & \vdots & \vdots  &
\g{3}  \textcolor{brown}{-\g{2N-1}}\\
\g{N}  \textcolor{brown}{-\g{N+2}}     &
\g{N-1}  \textcolor{brown}{-\g{N+1}}\textcolor{brown}{-\g{N+1}}                            &
\g{N-2}\textcolor{brown}{-\g{N+2}} \textcolor{brown}{-\g{N}}     &
\ldots & \ldots & \ldots &
2\g{1} \textcolor{brown}{-\g{2N-1}}\textcolor{brown}{-\g{3}}&
\g{0}+\g{2} \textcolor{brown}{-\g{2N}}\\
\hline
 \g{N+1}  \textcolor{brown}{-\g{N+1}}        &
 \g{N} \textcolor{brown}{-\g{N+2}} \textcolor{brown}{-\g{N}}     &
 \g{N-1} \textcolor{brown}{-\g{N+3}} \textcolor{brown}{-\g{N-1}}     &
 \ldots & \ldots  &
 \g{3}  \textcolor{brown}{-\g{2N-1}} \textcolor{brown}{-\g{3}}&
 \g{0}+\g{2} \textcolor{brown}{-\g{2N}}\textcolor{brown}{-\g{2}} \textcolor{red}{-\g{0}} &
 2\g{1} \textcolor{brown}{-\g{2N+1}}
\end{array}
\end{bmatrix}.
\]
We can rewrite $A_{0}^\textrm{anti}$ in several ways according to the type os structural properties we are looking for. For instance, }
$$
A_{0}^\textrm{anti} = S + B,
$$
where $S$ is a symmetric matrix and $B$ only has non-zero elements  on the first and last  rows and columns
in the following manner
$$
B  =
{\frac{1}{2}} \left[
\begin{array}{cccccc}
0 & \textcolor{red}{-g_0^\alpha} & 0& \dots & 0&0\\
g_2^\alpha & 0
& 0& \dots & 0 & g_{N}^\alpha \\
g_3^\alpha  & 0 & 0 & \dots &  & g_{N-1}^\alpha \\
\vdots & \vdots& \vdots &&\vdots&\vdots\\
g_{N-1}^\alpha  &0& 0 &&0& g_3^\alpha \\
g_N^\alpha  &0& 0 &&0& g_2^\alpha \\
0 & 0 & 0& \dots
&  \textcolor{red}{-g_{0}^\alpha}& 0
\end{array}
\right].
$$
%-------------------------------------------------------------------------------------------------------
{Alternatively, a different way of looking at the same matrix structure is put in evidence in the following equations.} We have
\begin{align*}
A_L^{\mathrm{anti}}  & = {T_L} - {\tilde{H}_L}  - \tilde{R}_L^{\mathrm{anti}}, \\%\label{eq:ALanti}\\
A_R^{\mathrm{anti}}  & = {T_R} - {\tilde{H}_R}  - \tilde{R}_R^{\mathrm{anti}}, %\label{eq:ARanti}
\end{align*}
where
\begin{equation} \label{eq:TLanti}
{T_L}  =
\begin{bmatrix} %\nonumber %\g{}
\g{1}   & \g{0}  & 0       & \ldots & 0 & 0 \\
\g{2}   & \g{1}  & \g{0}   & \ldots & 0 & 0 \\
\g{3}   & \g{2}  & \g{1}   & \ldots & 0 & 0 \\
\vdots  & \vdots & \vdots  & \vdots & \vdots & \vdots \\
\g{N+1} & \g{N}  & \g{N-1} & \ldots & \g{2} & \g{1}
\end{bmatrix},
\quad
{T_R}  = {T_L}^T
\end{equation}
are Toeplitz matrices in lower and upper Hessenberg form, respectively,
\begin{equation*}%\label{eq:HtildeLanti}
{\tilde{H}_L}  =
\begin{bmatrix} %\g{}
0   & \g{2}   & \g{3}   & \ldots & \g{N}   & \g{N+1} \\
0   & \g{3}   & \g{4}   & \ldots & \g{N+1} & \g{N+2} \\
0   & \g{4}   & \g{5}   & \ldots & \g{N+2} & \g{N+3} \\
\vdots  & \vdots & \vdots  & \vdots & \vdots & \vdots \\
0   & \g{N+2} & \g{N+3} & \ldots & \g{2N} & \g{2N+1}
\end{bmatrix},
\quad
{\tilde{H}_R}  = J{\tilde{H}_L}J,
\end{equation*}
$J$ being the antidiagonal or flip matrix, are Hankel matrices apart the first and last zero columns, respectively, and
\begin{equation*}%\label{eq:RtildeLanti}
\tilde{R}_L^{\mathrm{anti}}  =
\begin{bmatrix} %\g{}
0   & 0   & 0   & \ldots & 0 & 0 \\
0   & 0   & 0   & \ldots & 0 & 0 \\
0   & 0   & 0   & \ldots & 0 & 0 \\
\vdots  & \vdots & \vdots  & \vdots & \vdots & \vdots \\
0   & 0 & 0 & \ldots & \g{0} & 0
\end{bmatrix},
\quad
\tilde{R}_R^{\mathrm{anti}}  = J\tilde{R}_L^{\mathrm{anti}}J =
\begin{bmatrix} %\g{}
0   & \g{0}   & 0   & \ldots & 0 & 0 \\
0   & 0   & 0   & \ldots & 0 & 0 \\
0   & 0   & 0   & \ldots & 0 & 0 \\
\vdots  & \vdots & \vdots  & \vdots & \vdots & \vdots \\
0   & 0 & 0 & \ldots & 0 & 0
\end{bmatrix}
\end{equation*}
are just rank one correction matrices.\\
Thus, the matrix $A_0^{\mathrm{anti}}$ can we written as
\[
A_0^{\mathrm{anti}} =\frac{1}{2} \left( {T_L} + {T_R} -{{H}_L} -{{H}_R} \right) + R_0^{\mathrm{anti}}
\]
where ${{H}_L}$ and ${{H}_R}$ denote the full Hankel matrices linked to ${\tilde{H}_L}$ and ${\tilde{H}_L}$, respectively, and where
\[
R_0^{\mathrm{anti}} = \frac{1}{2}
\begin{bmatrix} %\g{}
0   & -\g{0}   & 0   & \ldots & 0 & \g{N+1} \\
\g{2}   & 0   & 0   & \ldots & 0 & \g{N} \\
\g{3}   & 0   & 0   & \ldots & 0 & \g{N-1} \\
\vdots  & \vdots & \vdots  & \vdots & \vdots & \vdots \\
\g{N}  & 0 & 0 & \ldots & 0 & \g{2}\\
\g{N+1}  & 0 & 0 & \ldots & -\g{0} & 0
\end{bmatrix}.
\]
Clearly, ${S_0}=\frac{1}{2} \left( {T_L} + {T_R} -{{H}_L} -{{H}_R} \right)$ is a symmetric matrix and, in particular, the matrix
\begin{equation*} %\label{eq:Toeplitz}
{T_0} =\frac{1}{2} \left ({T_L} + {T_R} \right )
\end{equation*}
is the real symmetric Toeplitz matrix with generating function $g_\alpha(\theta)$ (see \cite{DMS}  for the specific case and
\cite[Section 6.1]{GLT-BookI} for the general notion of generating function of Toeplitz matrices and Toeplitz matrix-sequences), whose Fourier coefficients are defined as
\begin{equation}
t_0=\g{1},\ t_1=(\g{0}+\g{2})/2,\ t_i=\g{i+1}/2, \ i=2, \ldots, N. \label{eq:Tcoeff}
\end{equation}
%
%}
From \cite{DMS}, we know that the generating function $g_\alpha(\theta)$ is nonnegative and not identically zero so that ${T_0}$ is positive definite for any matrix-size. Furthermore $g_\alpha(\theta)$ has a unique zero at $\theta=0$
and according to the order of the fractional derivative this zero has order $\alpha$ so that, by varying the matrix-sizes, the sequence of the minimal eigenvalues tends to zero as $m^{-\alpha}$ where $m=N+1$ is the matrix-size, in accordance with (\ref{eq:TLanti}) and taking into account \cite{S-extr2,BG-extr,S-extr1}.

%-------------------------------------------------------------------------------------------------------
%
\section{Anti-reflective numerical BCs}\label{sec:anti-refl}
A possible proposal to restore the continuity of the function and not only of its derivative is to consider numerical anti-reflective BCs, as
first introduced in the context of signal/image deblurring and  restoration \cite{S-AR-proposal}: see also \cite{imaging-AR-1,imaging-AR-2} for applications and further results in presence of noise and blurring. More precisely, we set
\begin{align}
u(-x+a,t)  - u(a,t) & = u(a,t) -u(x+a,t), \quad \mbox{for all} \quad x < a,
\label{Antireflective1}
\\
u(x+b,t) - u(b,t) & = u(b,t) -u(-x+b,t), \quad \mbox{for all} \quad x > b, \label{Antireflective2}
\end{align}
Thus, by considering the same arguments as before with respect to the boundaries, the approximation of the left fractional Riemann-Liouville derivative becomes
\begin{eqnarray}
 {}_{-\infty}^{RL}D_{x}^{\alpha,ref} u (x_j,t_n)
 &\approx & \frac{1}{(\Delta x)^\alpha} \left ( \sum_{k=0}^{j+1} g_{k}^\alpha U_{j+1-k}^n
 + \sum_{k=j+2}^{\infty} g_{k}^\alpha U_{j+1-k}^n \right )\nonumber \\
 & = & \frac{1}{(\Delta x)^\alpha} \left (\sum_{k=0}^{j+1} g_{k}^\alpha U_{j+1-k}^n
 + \sum_{k=j+2}^{\infty}  g_{k}^\alpha \left ( 2U_0^n - U_{k-(j+1)}^n \right ) \right )\nonumber  \\
 &\approx & \frac{1}{(\Delta x)^\alpha}\left ( \sum_{k=0}^{j+1} g_{k}^\alpha U_{j+1-k}^n
 + \sum_{k=j+2}^{N+j+1} g_{k}^\alpha \left ( 2U_0^n - U_{k-(j+1)}^n \right ) \right ) \label{glapproxl_antiR}
\end{eqnarray}
and the approximation of the right fractional Riemann-Liouville derivative becomes
\begin{eqnarray}
 {}_{x}^{RL}D_{\infty}^{\alpha,ref} u (x_j,t_n)
 &\approx & \frac{1}{(\Delta x)^\alpha} \left (\sum_{k=0}^{N-j+1} g_{k}^\alpha U_{j-1+k}^n
 + \sum_{k=N-j+2}^{\infty} g_{k}^\alpha  U_{j-1+k}^n \right ) \nonumber \\
 & = & \frac{1}{(\Delta x)^\alpha} \left (\sum_{k=0}^{N-j+1} g_{k}^\alpha U_{j-1+k}^n
 + \sum_{k=N-j+2}^{\infty} g_{k}^\alpha ( 2 U_N^n - U_{2N-j+1-k}^n) \right ) \nonumber \\
 & \approx  & \frac{1}{(\Delta x)^\alpha} \left ( \sum_{k=0}^{N-j+1} g_{k}^\alpha U_{j-1+k}^n
 + \sum_{k=N-j+2}^{2N-j+1} g_{k}^\alpha ( 2 U_N^n -U_{2N-j+1-k}^n) \right ). \label{glapproxr_antiR}
\end{eqnarray}
%-----------------------------------------------------------------------------------
%
Thus, by applying again the $\theta-$method, the matrix form of the problem is
\begin{equation*}
%{\bf U}^{n+1} = ( I +\mu_\alpha {\bf U}^{n},
(I -  \mu_\alpha \theta A_{\beta}^\textrm{antiR}) {\bf U}^{n+1} = (I +  \mu_\alpha (1-\theta) A_{\beta}^\textrm{antiR})  {\bf U}^{n},
\end{equation*}
with ${\bf U}^{n}=[ U_{0}^n, \dots, U_{N}^n]^T$, $I$ being the identity matrix, and $\mu_\alpha={\kappa_\alpha \Delta t}/{(\Delta x)^\alpha}$. By combining \eqref{glapproxl_antiR} and \eqref{glapproxr_antiR}, the matrix $A_{\beta}^\textrm{anti}$ is
expressed as the $(N+1)\times(N+1)$ matrix
\begin{equation}\label{eq:Abeta_antiR}
{A_{\beta}^\mathrm{antiR}} =  \frac{1+\beta}{2}{A_L^\mathrm{antiR}}+\frac{1-\beta}{2}{A_R^\mathrm{antiR}}
\end{equation}
{and for $\beta=0$ the matrix $2A_{0}^\textrm{antiR}$ equals the matrix
\[ \tiny \hskip -2cm
\begin{bmatrix}
\begin{array}{c|cccccc|c}
2\g{1} \textcolor{brown}{-\g{2N+1}} \textcolor{violet}{+2\g{0}+\z{1}} &
\g{0}+\g{2} \textcolor{brown}{-\g{2}} \textcolor{brown}{-\g{2N}} \textcolor{red}{-\g{0}} &
\g{3} \textcolor{brown}{-\g{3}} \textcolor{brown}{-\g{2N-1}} &
\ldots &  \ldots &
\g{N-1} \textcolor{brown}{-\g{N-1}} \textcolor{brown}{-\g{N+3}} &
\g{N} \textcolor{brown}{-\g{N}}\textcolor{brown}{-\g{N+2}}  &
\g{N+1} \textcolor{brown}{-\g{N+1}} \textcolor{violet}{+\z{N+1}}\\
\hline
\g{0}+\g{2}\textcolor{brown}{-\g{2N}} \textcolor{violet}{+\z{2}}&
2\g{1}\textcolor{brown}{-\g{3}} \textcolor{brown}{-\g{2N-1}}     &
\g{0}+\g{2} \textcolor{brown}{-\g{4}} \textcolor{brown}{-\g{2N-2}}&
\g{3}\textcolor{brown}{-\g{5}} \textcolor{brown}{-\g{2N-3}} &  \ldots      & \ldots &
\g{N-1} \textcolor{brown}{-\g{N+1}} \textcolor{brown}{-\g{N+1}}                     &
\g{N} \textcolor{brown}{-\g{N+2}} \textcolor{violet}{+\z{N}}    \\
\g{3}  \textcolor{brown}{-\g{2N-1}}  \textcolor{violet}{+\z{3}}   &
\g{0}+\g{2}\textcolor{brown}{-\g{4}} \textcolor{brown}{-\g{2N-2}}&
\vdots      & \vdots & \vdots & \vdots &
\g{N-2}\textcolor{brown}{-\g{N+2}}\textcolor{brown}{-\g{N}}                         &
\g{N-1} \textcolor{brown}{-\g{N+3}} \textcolor{violet}{+\z{N-1}}   \\
\vdots      & \vdots                               & \vdots      & \vdots & \vdots & \vdots & \vdots                          & \vdots      \\
\g{N-1} \textcolor{brown}{-\g{N+3}} \textcolor{violet}{+\z{N-1}}   &
\g{N-2}\textcolor{brown}{-\g{N}} \textcolor{brown}{-\g{N+2}}    &
\vdots                          & \vdots & \vdots & \vdots & \vdots  &
\g{3}  \textcolor{brown}{-\g{2N-1}} \textcolor{violet}{+\z{3}}\\
\g{N}  \textcolor{brown}{-\g{N+2}}     \textcolor{violet}{+\z{N}}&
\g{N-1}  \textcolor{brown}{-\g{N+1}}\textcolor{brown}{-\g{N+1}}                            &
\g{N-2}\textcolor{brown}{-\g{N+2}} \textcolor{brown}{-\g{N}}     &
\ldots & \ldots & \ldots &
2\g{1} \textcolor{brown}{-\g{2N-1}}\textcolor{brown}{-\g{3}}&
\g{0}+\g{2} \textcolor{brown}{-\g{2N}}\textcolor{violet}{+\z{2}}\\
\hline
 \g{N+1}  \textcolor{brown}{-\g{N+1}}    \textcolor{violet}{+\z{N+1}}    &
 \g{N} \textcolor{brown}{-\g{N+2}} \textcolor{brown}{-\g{N}}     &
 \g{N-1} \textcolor{brown}{-\g{N+3}} \textcolor{brown}{-\g{N-1}}     &
 \ldots & \ldots  &
 \g{3}  \textcolor{brown}{-\g{2N-1}} \textcolor{brown}{-\g{3}}&
 \g{0}+\g{2} \textcolor{brown}{-\g{2N}}\textcolor{brown}{-\g{2}} \textcolor{red}{-\g{0}} &
 2\g{1} \textcolor{brown}{-\g{2N+1}} \textcolor{violet}{+2\g{0}+\z{1}}
\end{array}
\end{bmatrix},
\]
with
\[
\z{r} =2 \sum_{k=r+1}^{N+r} \g{k}, \quad r=1, \ldots, N+1.
\]
}
%--------------------------------------------C
More in detail the obtained structure takes the more explicit form
\begin{align*}
A_L^{\mathrm{antiR}}  & = {T_L} -{\tilde{H}_L}  - \hat{R}_L^{\mathrm{anti}}+ \hat{Z}_L^{\mathrm{antiR}}, \\% \label{eq:ALantiR}\\
A_R^{\mathrm{antiR}}  & = {T_R} - {\tilde{H}_R}  - \hat{R}_R^{\mathrm{anti}}+ \hat{Z}_L^{\mathrm{anti}R}, %\label{eq:ARantiR}
\end{align*}
where
\begin{equation*}%\label{eq:RtildeantiR}
\hat{R}_L^{\mathrm{antiR}}  =
\begin{bmatrix} %\g{}
0   & 0   & 0   & \ldots & 0 & 0 \\
0   & 0   & 0   & \ldots & 0 & 0 \\
0   & 0   & 0   & \ldots & 0 & 0 \\
\vdots  & \vdots & \vdots  & \vdots & \vdots & \vdots \\
0   & 0 & 0 & \ldots & \g{0} & -2\g{0}
\end{bmatrix},
\
\hat{R}_R^{\mathrm{antiR}}  = J\tilde{R}_L^{\mathrm{antiR}}J
\end{equation*}
and
\begin{equation*}%\label{eq:ZtildeantiR}
\hat{Z}_L^{\mathrm{antiR}}  =
\begin{bmatrix} %\g{}
\ 2{\tilde{H}_L} \mathbf{e}\ | O^{(N+1) \times N}
\end{bmatrix},
\quad
\hat{Z}_R^{\mathrm{antiR}}  =
\begin{bmatrix} %\g{}
O^{(N+1) \times N} \ | \ 2{\tilde{H}_R} \mathbf{e}
\end{bmatrix},
\end{equation*}
with $\mathbf{e}=[1, \ldots, 1]^T$ and $O^{(N+1) \times N}$ zero matrix of dimension $(N+1)\times N$. \\
Thus, as in the case of numerical anti-symmetric BCs, for $\beta=0$ the matrix $A_0^{\mathrm{antiR}}$ can we written as
\[
A_0^{\mathrm{antiR}}  =\frac{1}{2} \left( {T_L} + {T_R} -{{H}_L} -{{H}_R} \right) + R_0^{\mathrm{antiR}} = {S_0} + R_0^{\mathrm{antiR}}\\
\]
where again ${{H}_L}$ and ${{H}_R}$ denote the full Hankel matrices linked to ${\tilde{H}_L}$ and ${\tilde{H}_L}$ respectively and where
\[
R_0^{\mathrm{antiR}} = \frac{1}{2}
\begin{bmatrix} %\g{}
\begin{array}{c|cccc|c}
2\g{0} +\z{1}   & -\g{0}   & 0   & \ldots & 0 & \g{N+1} + \z{N+1} \\
\hline
\g{2} + \z{2} & 0   & 0   & \ldots & 0 & \g{N} + \z{N}\\
\g{3} + \z{3}  & 0   & 0   & \ldots & 0 & \g{N-1} \\
\vdots  & \vdots & \vdots  & \vdots & \vdots & \vdots \\
\g{N} + \z{N} & 0 & 0 & \ldots & 0 & \g{2}+ \z{2}\\
\hline
\g{N+1}+ \z{N+1}  & 0 & 0 & \ldots & -\g{0} & 2\g{0}+ \z{1}
\end{array}
\end{bmatrix}.
\]
%--------------------------------------------------------------------------------------------------------
%
\section{Numerical experiments}\label{sec:num}
\begin{figure}[h]
{\center
\includegraphics[width=\textwidth]{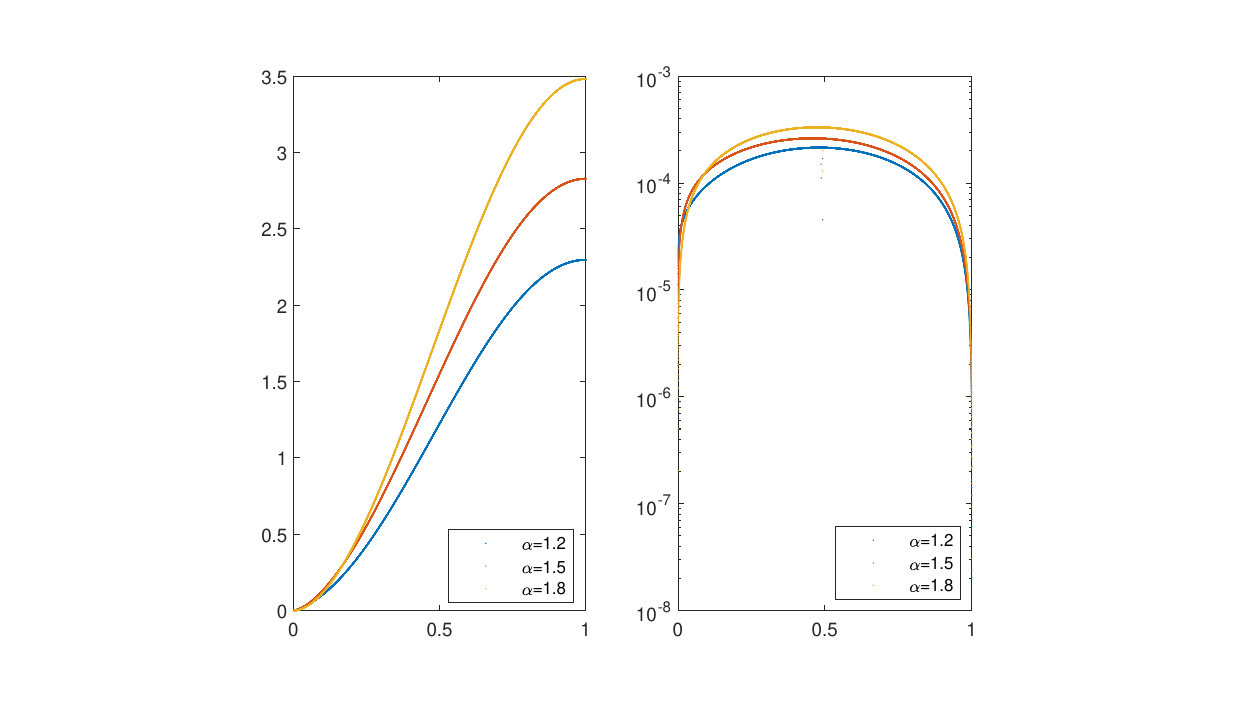}
}
\caption{Eigenvalue distribution of ${A_{0}^\mathrm{anti}}$ with size $16000$ for $\alpha=1.2, 1.5, 1.8$ and absolute error with respect to the generating function $f_{\alpha,T_{0}}(\theta)$.}
\label{fig:distribution}
\end{figure}

In the following, in view of better approximation/stability properties, we consider numerical experiments related the more relevant case of application of Implicit Euler scheme or Crank-Nicolson scheme
to the solution of \eqref{eq:pb}-\eqref{reflective2}. Indeed, in such case the matrix form of the problem is
\begin{equation} \label{eq:EI}
\left ({I} - \mu_\alpha{A_{\beta}} \right ) {\bf U}^{n+1} = {\bf U}^{n},
\end{equation}
and
\begin{equation} \label{eq:CN}
\left ({I} - \frac{\mu_\alpha}{2}{A_{\beta}} \right ) {\bf U}^{n+1} =
                \left ( {I} + \frac{\mu_\alpha}{2} {A_{\beta}} \right ) {\bf U}^{n},
\end{equation}
respectively, where ${A_{\beta}}$ is as in \eqref{eq:Abeta_anti} or as in \eqref{eq:Abeta_antiR}. \\
Since the implicit schemes involve the solution of a linear system, in the following we will focus on the application of Krylov methods and related effective preconditioned techniques in the case $\beta=0$. \\
On that point of view,  we start by giving numerical evidence of the spectral analysis of the involved structured matrices in the case $\beta=0$,
namely {${T_{0}}$, ${A_{0}^\mathrm{anti}}$, and ${A_{0}^\mathrm{antiR}}$}.
First of all, we highlight as the minimal eigenvalues goes to zero asymptotically as $m^{-\alpha}$, $m$ being the matrix dimension, according to the order of zero of the Toeplitz generating function (see \cite{DMS} and references therein)
\[
f_{\alpha,{T_{0}}}(\theta)= - (f_{\alpha,{T_L}}(\theta) + f_{\alpha,{T_R}}(\theta))/2
\]
where
\begin{align*}
f_{\alpha,{T_L}}(\theta) &= \sum_{k=-1}^\infty \g{k+1} e^{ik\theta} = {e^{-i\theta}}(1+ e^{i(\theta+\pi)})^\alpha \\
f_{\alpha,{T_R}}(\theta) &= \overline{f_{\alpha,{T_L}}(\theta)}. \\
\end{align*}
Indeed  in Table \ref{tab:eigmin} we report the minimal eigenvalue of the matrix $X_m\in \mathbb{R}^{m\times m}$, with {$X \in { {T_{0}},  {A_{0}^\mathrm{anti}}, {A_{0}^\mathrm{antiR}} }$}  for increasing dimension together
with the corresponding quantity
\[
\gamma(X_m)= \log_2\left ( \frac{\lambda_{\min} (X_m)} {\lambda_{\min} (X_{2m})} \right ),
\]
whose limit for the dimension $m$ tending to infinity is exactly $\alpha$. In fact, the latter is expected since the related Toeplitz matrix admits a generating function in the Wiener class which is nonnegative and with a unique zero of order $\alpha$ at $\theta=0$. Hence, in the light of the results in \cite{S-extr1,S-extr2,BG-extr}, we know that the minimal eigenvalue of $T_m(f_{\alpha,{T_{0}}})$ is asymptotic to $m^{-\alpha}$. We notice that the result in \cite{BG-extr} is far more general, where the result concerns the minimal modulus of the eigenvalues and where the assumption of a nonnegative generating function is replaced by a complex-valued generating function with weak sectorial character.\\
In addition, Figure \ref{fig:distribution}.a highlights how the whole eigenvalue distribution of the matrix ${A_{0}^\mathrm{anti}}$ mimics in quite good measure the quoted generating function even in the case of a moderate matrix dimension as $16000$, while for the sake of completeness in Figure \ref{fig:distribution}.b the absolute error with respect the generating function is plotted for different values of the parameter $\alpha$.
Regarding the absolute error reported in Figure \ref{fig:distribution}.b, it is interesting to observe the smooth shape of the error which would suggest the use of the extrapolation techniques for the fast eigenvalue computation as in \cite{extrap1,extrap2}: we also notice the exception of a unique double eigenvalue with much smaller error close to the abscissa value of $0.5$.
\ \\
Then, we consider the spectral analysis of the whole matrix $\mathcal{A}_0^{\mathrm{anti}}={I} - \mu_\alpha{A_{0}^\mathrm{anti}}$ and of the
 Toeplitz counterpart ${\mathcal{T}_0}={I} - \mu_\alpha {T_{0}}$.
In Table \ref{tab:EI_lAlT}  the minimal and maximal eigenvalues are reported, together with the spectral condition number $K_2$ for increasing dimensions
in the case $k=1$ and $\Delta t=\Delta x$ and Implicit Euler scheme  for different values of the parameter $\alpha$.
The same analysis is considered in Table \ref{tab:CN_lAlT} with respect to the Crank-Nicolson scheme. \\
It is evident the worsening of the condition number for increasing dimension as the parameter $\alpha$ increases, approaching the standard second order differential problem case, due to a faster decrease of the minimal eigenvalue to zero, while the maximal one tends to a constant (the maximum of the generating function which is known to be continuous and $2\pi$-periodic \cite{DMS}). \\
Then, we consider the GMRES method alone or with suitable preconditioners taken in the algebra of the circulant matrices and in
the algebra of sine transforms (also called $\tau$ matrices \cite{BC}) for the solution of a linear system with matrix $\mathcal{A}_0^{\mathrm{anti}}$. More precisely, we compare the case of no preconditioning, Strang Circulant preconditioning, Frobenius Optimal Circulant preconditioning, natural $\tau$ preconditioning, Frobenius Optimal $\tau$ preconditioning (see \cite{CN,DiB-S,sh,tau-theory2} and references therein),
{all built by referring to the Toeplitz part $T_0$.}\\
In Table \ref{tab:EICN_gmresit} the number of iterations required to reach convergence within a tolerance of $10^{-6}$ is reported in the case $k=1$ and $\Delta t=\Delta x$ and both Implicit Euler and Crank-Nicolson schemes for different values of the parameter $\alpha$.
The constant number of iterations independent of the dimension testify the effectiveness and robustness of the proposed preconditioning strategies. An essentially negligible dependence on the chosen scheme is also observed. \\
Finally, we test the robustness also when increasing the value $k$.
In Table \ref{tab:EICN_gmresit_k100} the very same analysis is reported, giving evidence of the strong robustness of Strang Circulant and $\tau$ preconditioners.
In Tables \ref{tab:EICN_lAlT_AntiR}, \ref{tab:EIcn_gmresit_AntiR}, and \ref{tab:EIcn_gmresit_AntiR_k100} we collect the same tests in the case of
numerical anti-reflective BCs and the comments are very similar.
Finally, Table \ref{tab:EICN_gmresit_random} and Table \ref{tab:EIcn_gmresit_AntiR_random} account for the use of a random solution for numerical anti-symmetric and anti-reflective BCs, respectively:  the randomness is nonphysical and is introduced for checking the numerical robustness of the considered methods and the conclusion is that no changes are observed in the related convergence history.
%------------------------------------------------
\begin{table}
  \centering
  \footnotesize
  \begin{tabular}{|c|cc|cc|cc|}
    \hline
   $m$ & $\lambda_{\min}({T_0})$ & $\gamma({T_0})$
    & $\lambda_{\min}({A_0^{\mathrm{anti}}})$ & $\gamma({A_0^{\mathrm{anti}}})$
    & $\lambda_{\min}({A_0^{\mathrm{antiR}}})$ & $\gamma({A_0^{\mathrm{antiR}}})$\\
    \hline
%----------------------------------------------------------------------------------------------
    \multicolumn {7}{|c|}{$\alpha=1.2$} \\
    \hline
1000 & 2.33167e-04 & -           & 3.15528e-04 & -           & 3.11001e-05 & - \\
2000 & 1.01115e-04 & 1.20536 & 1.37002e-04 & 1.20358 & 1.35354e-05 & 1.20019 \\
4000 & 4.39242e-05 & 1.20292 & 5.95594e-05 & 1.20179 & 5.89123e-06 & 1.20009 \\
8000 & 1.90983e-05 & 1.20158 & 2.59087e-05 & 1.20090 & 2.56422e-06 & 1.20005 \\
\hline
%----------------------------------------------------------------------------------------------
    \multicolumn {7}{|c|}{$\alpha=1.5$} \\
    \hline
1000 & 1.01144e-04 & -       & 1.26438e-04 & -       & 6.05846e-06 & \\
2000 & 3.57435e-05 & 1.50066 & 4.46521e-05 & 1.50164 & 2.14144e-06 & 1.50037 \\
4000 & 1.26338e-05 & 1.50039 & 1.57779e-05 & 1.50082 & 7.57015e-07 & 1.50019 \\
8000 & 4.46604e-06 & 1.50023 & 5.57676e-06 & 1.50041 & 2.67628e-07 & 1.50009 \\
   \hline
%----------------------------------------------------------------------------------------------
    \multicolumn {7}{|c|}{$\alpha=1.8$} \\
    \hline
    1000 & 2.69766e-05 &      -  & 2.99102e-05 & -       & 4.47444e-07 & -  \\
    2000 & 7.75208e-06 & 1.79905 & 8.58130e-06 & 1.80137 & 1.28440e-07 & 1.80061 \\
    4000 & 2.22692e-06 & 1.79954 & 2.46316e-06 & 1.80069 & 3.68771e-08 & 1.80030 \\
    8000 & 6.39615e-07 & 1.79977 & 7.07189e-07 & 1.80034 & 1.05890e-08 & 1.80015 \\
   \hline
%----------------------------------------------------------------------------------------------
  \end{tabular}
  \caption{Spectral Analysis of
 {${T_0}$, $A_0^{\mathrm{anti}}$, and
  $A_0^{\mathrm{antiR}}$}    - Implicit Euler method case - case $k=1$.} \label{tab:eigmin}
\end{table}
%------------------------------------------------
%------------------------------------------------
\begin{table}
  \centering
  \footnotesize
  \begin{tabular}{|c|ccc|ccc|}
    \hline
    $n$ &
    $\lambda_{\min}(\mathcal{A}_0^{\mathrm{anti}})$ & $\lambda_{\max}(\mathcal{A}_0^{\mathrm{anti}})$ & $K_2(\mathcal{A}_0^{\mathrm{anti}})$ &
          $\lambda_{\min}({\mathcal{T}_0})$ & $\lambda_{\max}({\mathcal{T}_0})$ & $K_2({\mathcal{T}_0})$ \\
    \hline
%----------------------------------------------------------------------------------------------
    \multicolumn {7}{|c|}{$\alpha=1.2$} \\
    \hline
1000  & 1.00125e+00 &  1.01443e+01 & 1.01315e+01 & 1.00093e+00 & 1.01443e+01 & 1.01349e+01  \\
2000  & 1.00063e+00 &  1.15051e+01 & 1.14979e+01 & 1.00046e+00 & 1.15051e+01 & 1.14997e+01  \\
4000  & 1.00031e+00 &  1.30677e+01 & 1.30637e+01 & 1.00023e+00 & 1.30677e+01 & 1.30647e+01  \\
8000  & 1.00016e+00 &  1.48625e+01 & 1.48602e+01 & 1.00012e+00 & 1.48625e+01 & 1.48608e+01  \\
\hline
%----------------------------------------------------------------------------------------------
    \multicolumn {7}{|c|}{$\alpha=1.5$} \\
    \hline
1000 & 1.00399e+00 &  9.03978e+01 & 9.00385e+01 & 1.00319e+00 & 9.03978e+01 & 9.01102e+01  \\
2000 & 1.00199e+00 &  1.27459e+02 & 1.27206e+02 & 1.00160e+00 & 1.27459e+02 & 1.27256e+02  \\
4000 & 1.00100e+00 &  1.79863e+02 & 1.79684e+02 & 1.00080e+00 & 1.79863e+02 & 1.79720e+02  \\
8000 & 1.00050e+00 &  2.53966e+02 & 2.53840e+02 & 1.00040e+00 & 2.53966e+02 & 2.53865e+02  \\
   \hline
%----------------------------------------------------------------------------------------------
    \multicolumn {7}{|c|}{$\alpha=1.8$} \\
    \hline
1000 & 1.00749e+00 &  8.74988e+02 & 8.68480e+02 & 1.00676e+00 & 8.74988e+02 & 8.69114e+02  \\
2000 & 1.00375e+00 &  1.52331e+03 & 1.51762e+03 & 1.00339e+00 & 1.52331e+03 & 1.51817e+03  \\
4000 & 1.00187e+00 &  2.65203e+03 & 2.64707e+03 & 1.00169e+00 & 2.65203e+03 & 2.64755e+03  \\
8000 & 1.00094e+00 &  4.61718e+03 & 4.61285e+03 & 1.00085e+00 & 4.61718e+03 & 4.61327e+03  \\
   \hline
%----------------------------------------------------------------------------------------------
  \end{tabular}
  \caption{Spectral Analysis of $\mathcal{A}_0^{\mathrm{anti}}= I-\mu_\alpha A_0^{\mathrm{anti}}$ and ${\mathcal{T}_0}=I-\mu_\alpha {T_0}$ - Implicit Euler method case - case $k=1$.} \label{tab:EI_lAlT}
\end{table}
\begin{table}
  \centering
  \footnotesize
  \begin{tabular}{|c|ccc|ccc|}
    \hline
    $n$ &
    $\lambda_{\min}(\mathcal{A}_0^{\mathrm{anti}})$ & $\lambda_{\max}(\mathcal{A}_0^{\mathrm{anti}})$ & $K_2(\mathcal{A}_0^{\mathrm{anti}})$ &
          $\lambda_{\min}({\mathcal{T}_0})$ & $\lambda_{\max}({\mathcal{T}_0})$ & $K_2({\mathcal{T}_0})$ \\
    \hline
%----------------------------------------------------------------------------------------------
    \multicolumn {7}{|c|}{$\alpha=1.2$} \\
    \hline
1000 & 1.00063e+00 &  5.57213e+00 & 5.56863e+00 & 1.00046e+00 & 5.57213e+00 & 5.56954e+00  \\
2000 & 1.00031e+00 &  6.25253e+00 & 6.25057e+00 & 1.00023e+00 & 6.25253e+00 & 6.25108e+00  \\
4000 & 1.00016e+00 &  7.03387e+00 & 7.03277e+00 & 1.00012e+00 & 7.03387e+00 & 7.03306e+00  \\
8000 & 1.00008e+00 &  7.93127e+00 & 7.93066e+00 & 1.00006e+00 & 7.93127e+00 & 7.93082e+00  \\
\hline
%----------------------------------------------------------------------------------------------
    \multicolumn {7}{|c|}{$\alpha=1.5$} \\
    \hline
1000 & 1.00200e+00 &  4.56989e+01 & 4.56079e+01 & 1.00160e+00 & 4.56989e+01 & 4.56261e+01  \\
2000 & 1.00100e+00 &  6.42297e+01 & 6.41657e+01 & 1.00080e+00 & 6.42297e+01 & 6.41785e+01  \\
4000 & 1.00050e+00 &  9.04315e+01 & 9.03865e+01 & 1.00040e+00 & 9.04315e+01 & 9.03954e+01  \\
8000 & 1.00025e+00 &  1.27483e+02 & 1.27451e+02 & 1.00020e+00 & 1.27483e+02 & 1.27458e+02  \\
   \hline
%----------------------------------------------------------------------------------------------
    \multicolumn {7}{|c|}{$\alpha=1.8$} \\
    \hline
1000 & 1.00375e+00 &  4.37994e+02 & 4.36359e+02 & 1.00338e+00 & 4.37994e+02 & 4.36519e+02  \\
2000 & 1.00187e+00 &  7.62157e+02 & 7.60731e+02 & 1.00169e+00 & 7.62157e+02 & 7.60868e+02  \\
4000 & 1.00094e+00 &  1.32652e+03 & 1.32527e+03 & 1.00085e+00 & 1.32652e+03 & 1.32539e+03  \\
8000 & 1.00047e+00 &  2.30909e+03 & 2.30801e+03 & 1.00042e+00 & 2.30909e+03 & 2.30811e+03  \\
   \hline
%----------------------------------------------------------------------------------------------
  \end{tabular}
  \caption{Spectral Analysis of $\mathcal{A}_0^{\mathrm{anti}}= I-\frac{\mu_\alpha}{2} A_0^{\mathrm{anti}}$ and ${\mathcal{T}_0}=I-\frac{\mu_\alpha}{2} {T_0}$ - Crank-Nicolson method case - case $k=1$.} \label{tab:CN_lAlT}
\end{table}
\begin{table}
  \centering
  \footnotesize
  \begin{tabular}{|c|ccccc|ccccc|}
  \hline
    &  \multicolumn {5}{c|}{Implicit Euler} &  \multicolumn {5}{c|}{Crank-Nicolson}\\
    \hline
    $n$
    & - & $\mathcal{C}_0$ & $\mathcal{C}_{0,\mathrm{ott}}$ & $\mathcal{\tau}_0$ & $\mathcal{\tau}_{0,\mathrm{ott}}$
    & - & $\mathcal{C}_0$ & $\mathcal{C}_{0,\mathrm{ott}}$ & $\mathcal{\tau}_0$ & $\mathcal{\tau}_{0,\mathrm{ott}}$ \\
    \hline
%----------------------------------------------------------------------------------------------
    \multicolumn {11}{|c|}{$\alpha=1.2$} \\
    \hline
1000 & 21 &  5 & 5 & 3 &  3  & 15 &  5 & 5 & 3 &  3  \\
2000 & 22 &  5 & 5 & 3 &  3  & 15 &  5 & 5 & 3 &  3  \\
4000 & 23 &  5 & 5 & 3 &  3  & 16 &  5 & 5 & 3 &  3  \\
8000 & 24 &  5 & 5 & 3 &  3  & 17 &  5 & 5 & 3 &  3  \\
   \hline
%----------------------------------------------------------------------------------------------
    \multicolumn {11}{|c|}{$\alpha=1.5$} \\
    \hline
1000 & 63 &  5 & 5 & 4 &  4  & 46 &  5 & 5 & 4 &  4  \\
2000 & 74 &  5 & 5 & 4 &  4  & 53 &  5 & 5 & 4 &  4  \\
4000 & 87 &  5 & 5 & 4 &  4  & 62 &  5 & 5 & 4 &  4  \\
8000 & 101 &  5 & 5 & 4 &  4 & 73 &  5 & 5 & 4 &  4  \\
     \hline
%----------------------------------------------------------------------------------------------
    \multicolumn {11}{|c|}{$\alpha=1.8$} \\
    \hline
1000 & 172 &  5 & 7 & 4 &  4  & 127 &  5 & 6 & 4 &  4  \\
2000 & 220 &  5 & 7 & 4 &  4  & 163 &  5 & 6 & 3 &  4  \\
4000 & 281 &  5 & 6 & 3 &  4  & 208 &  5 & 6 & 3 &  4  \\
8000 & 358 &  5 & 6 & 3 &  4  & 265 &  5 & 5 & 3 &  3  \\
     \hline
  \end{tabular}
  \caption{Number of preconditioned GMRES iterations to solve the linear system with matrix $\mathcal{A}_0^{\mathrm{anti}}=  \nu I-{k}A_0^{\mathrm{anti}}$  (Implicit Euler method) and $\mathcal{A}_0^{\mathrm{anti}}=  \nu I-\frac{k}{2}A_0^{\mathrm{anti}}$ (Crank-Nicolson method) for increasing dimension $n$  till $tol =1.e-6$ - case $k=1$.} \label{tab:EICN_gmresit}
\end{table}
\begin{table}
  \centering
  \footnotesize
  \begin{tabular}{|c|ccccc|ccccc|}
  \hline
    &  \multicolumn {5}{c|}{Implicit Euler} &  \multicolumn {5}{c|}{Crank-Nicolson}\\
    \hline
    $n$
    & - & $\mathcal{C}_0$ & $\mathcal{C}_{0,\mathrm{ott}}$ & $\mathcal{\tau}_0$ & $\mathcal{\tau}_{0,\mathrm{ott}}$
    & - & $\mathcal{C}_0$ & $\mathcal{C}_{0,\mathrm{ott}}$ & $\mathcal{\tau}_0$ & $\mathcal{\tau}_{0,\mathrm{ott}}$ \\
    \hline
%----------------------------------------------------------------------------------------------
    \multicolumn {11}{|c|}{$\alpha=1.2$} \\
    \hline
1000 & 167 &  6 & 8 & 4 &  4  & 128 &  6 & 7 & 4 &  4  \\
2000 & 187 &  6 & 7 & 4 &  4  & 141 &  6 & 6 & 4 &  4  \\
4000 & 205 &  6 & 7 & 4 &  4  & 152 &  6 & 6 & 4 &  4  \\
8000 & 227 &  6 & 6 & 4 &  4  & 163 &  6 & 6 & 4 &  4  \\
   \hline
%----------------------------------------------------------------------------------------------
    \multicolumn {11}{|c|}{$\alpha=1.5$} \\
    \hline
1000 & 359 &  5 & 13 & 4 &  4  & 306 &  5 & 11 & 4 &  4  \\
2000 & 481 &  5 & 13 & 4 &  4  & 381 &  5 & 11 & 4 &  4  \\
4000 & 599 &  5 & 12 & 4 &  4  & 472 &  5 & 10 & 4 &  4  \\
8000 & 745 &  5 & 11 & 4 &  4  & 566 &  5 & 9 & 4 &  4  \\
     \hline
%----------------------------------------------------------------------------------------------
    \multicolumn {11}{|c|}{$\alpha=1.8$} \\
    \hline
1000 & 705  &  4 & 20 & 4 &  5  & 666  &  4 & 19 & 4 &  4  \\
2000 & 1158 &  4 & 24 & 4 &  4  & 974  &  4 & 22 & 4 &  4  \\
4000 & 1597 &  4 & 27 & 4 &  4  & 1303 &  4 & 22 & 4 &  4  \\
8000 & 2249 &  4 & 27 & 4 &  4  & 1787 &  4 & 21 & 3 &  4  \\
     \hline
  \end{tabular}
  \caption{Number of preconditioned GMRES iterations to solve the linear system with matrix $\mathcal{A}_0^{\mathrm{anti}}=  \nu I-{k}A_0^{\mathrm{anti}}$  (Implicit Euler method) and $\mathcal{A}_0^{\mathrm{anti}}=  \nu I-\frac{k}{2}A_0^{\mathrm{anti}}$ (Crank-Nicolson method) for increasing dimension $n$  till $tol =1.e-6$ - case $k=100$.} \label{tab:EICN_gmresit_k100}
\end{table}
\begin{table}
  \centering
  \footnotesize
  \begin{tabular}{|c|ccccc|ccccc|}
  \hline
    &  \multicolumn {5}{c|}{Implicit Euler} &  \multicolumn {5}{c|}{Crank-Nicolson}\\
    \hline
    $n$
    & - & $\mathcal{C}_0$ & $\mathcal{C}_{0,\mathrm{ott}}$ & $\mathcal{\tau}_0$ & $\mathcal{\tau}_{0,\mathrm{ott}}$
    & - & $\mathcal{C}_0$ & $\mathcal{C}_{0,\mathrm{ott}}$ & $\mathcal{\tau}_0$ & $\mathcal{\tau}_{0,\mathrm{ott}}$ \\
    \hline
%----------------------------------------------------------------------------------------------
    \multicolumn {11}{|c|}{$\alpha=1.2$} \\
    \hline
1000 & 21 &  5 & 5 & 3 &  3 & 16 &  5 & 5 & 3 &  3  \\
2000 & 22 &  5 & 5 & 3 &  3 & 16 &  5 & 5 & 3 &  3  \\
4000 & 23 &  5 & 5 & 3 &  3 & 17 &  5 & 5 & 3 &  3  \\
8000 & 24 &  5 & 5 & 3 &  3 & 18 &  5 & 5 & 3 &  3  \\
   \hline
%----------------------------------------------------------------------------------------------
    \multicolumn {11}{|c|}{$\alpha=1.5$} \\
    \hline
1000 & 54 &  5 & 5 & 4 &  4  & 41 &  5 & 5 & 4 &  4  \\
2000 & 62 &  5 & 5 & 4 &  4  & 46 &  5 & 5 & 4 &  4  \\
4000 & 70 &  5 & 5 & 4 &  4  & 53 &  5 & 5 & 4 &  4  \\
8000 & 80 &  5 & 5 & 4 &  4  & 61 &  5 & 5 & 4 &  4  \\
     \hline
%----------------------------------------------------------------------------------------------
    \multicolumn {11}{|c|}{$\alpha=1.8$} \\
    \hline
1000 & 136 &  4 & 7 & 4 &  4  & 104 &  5 & 6 & 4 &  4  \\
2000 & 166 &  5 & 7 & 4 &  4  & 126 &  5 & 6 & 4 &  4  \\
4000 & 200 &  5 & 6 & 3 &  4  & 152 &  5 & 6 & 3 &  4  \\
8000 & 242 &  5 & 6 & 4 &  4  & 187 &  5 & 6 & 4 &  4  \\
     \hline
  \end{tabular}
  \caption{Number of preconditioned GMRES iterations to solve the linear system with matrix $\mathcal{A}_0^{\mathrm{anti}}=  \nu I-{k}A_0^{\mathrm{anti}}$  (Implicit Euler method) and $\mathcal{A}_0^{\mathrm{anti}}=  \nu I-\frac{k}{2}A_0^{\mathrm{anti}}$ (Crank-Nicolson method) for increasing dimension $n$  till $tol =1.e-6$ - case $k=1$ - random exact solution.} \label{tab:EICN_gmresit_random}
\end{table}
%%
%--------------------------------------ANTIR-----------------------------------------------------------------
%
\begin{table}
  \centering
  \footnotesize
  \begin{tabular}{|c|ccc|ccc|}
    \hline
      &  \multicolumn {3}{c|}{Implicit Euler} &  \multicolumn {3}{c|}{Crank-Nicolson}\\
    \hline
    $n$
    & $\lambda_{\min}(\mathcal{A}_0^{\mathrm{antiR}})$ & $\lambda_{\max}(\mathcal{A}_0^{\mathrm{antiR}})$ & $K_2(\mathcal{A}_0^{\mathrm{antiR}})$
    & $\lambda_{\min}(\mathcal{A}_0^{\mathrm{antiR}})$ & $\lambda_{\max}(\mathcal{A}_0^{\mathrm{antiR}})$ & $K_2(\mathcal{A}_0^{\mathrm{antiR}})$  \\
    \hline
%----------------------------------------------------------------------------------------------
    \multicolumn {7}{|c|}{$\alpha=1.2$} \\
    \hline
1000 & 1.00012e+00 &  1.01443e+01 & 1.33754e+01 & 1.00006e+00 & 5.57213e+00 & 6.61947e+00  \\
2000 & 1.00006e+00 &  1.15051e+01 & 1.55581e+01 & 1.00003e+00 & 6.25253e+00 & 7.56384e+00  \\
4000 & 1.00003e+00 &  1.30677e+01 & 1.81494e+01 & 1.00002e+00 & 7.03387e+00 & 8.67653e+00  \\
8000 & 1.00002e+00 &  1.48625e+01 & 2.12311e+01 & 1.00001e+00 & 7.93127e+00 & 9.98983e+00  \\
\hline
%----------------------------------------------------------------------------------------------
    \multicolumn {7}{|c|}{$\alpha=1.5$} \\
    \hline
1000 & 1.00019e+00 &  9.03978e+01 & 1.99489e+02 & 1.00010e+00 & 4.56989e+01 & 8.38871e+01  \\
2000 & 1.00010e+00 &  1.27459e+02 & 3.10333e+02 & 1.00005e+00 & 6.42297e+01 & 1.29009e+02  \\
4000 & 1.00005e+00 &  1.79863e+02 & 4.84674e+02 & 1.00002e+00 & 9.04315e+01 & 1.99620e+02  \\
8000 & 1.00002e+00 &  2.53966e+02 & 7.59411e+02 & 1.00001e+00 & 1.27483e+02 & 3.10438e+02  \\
   \hline
%----------------------------------------------------------------------------------------------
    \multicolumn {7}{|c|}{$\alpha=1.8$} \\
    \hline
1000 & 1.00011e+00 &  8.74988e+02 & 2.92340e+03 & 1.00006e+00 & 4.37994e+02 & 1.22410e+03  \\
2000 & 1.00006e+00 &  1.52331e+03 & 5.89058e+03 & 1.00003e+00 & 7.62157e+02 & 2.45650e+03  \\
4000 & 1.00003e+00 &  2.65203e+03 & 1.18932e+04 & 1.00001e+00 & 1.32652e+03 & 4.94501e+03  \\
8000 & 1.00001e+00 &  4.61718e+03 & 2.40495e+04 & 1.00001e+00 & 2.30909e+03 & 9.97769e+03  \\
   \hline
%----------------------------------------------------------------------------------------------
  \end{tabular}
  \caption{Spectral Analysis of $\mathcal{A}_0^{\mathrm{antiR}}$  - Implicit Euler and Crank-Nicolson method case - case $k=1$.} \label{tab:EICN_lAlT_AntiR}
\end{table}
\begin{table}
  \centering
  \footnotesize
  \begin{tabular}{|c|ccccc|ccccc|}
  \hline
    &  \multicolumn {5}{c|}{Implicit Euler} &  \multicolumn {5}{c|}{Crank-Nicolson}\\
    \hline
    $n$
    & - & $\mathcal{C}_0$ & $\mathcal{C}_{0,\mathrm{ott}}$ & $\mathcal{\tau}_0$ & $\mathcal{\tau}_{0,\mathrm{ott}}$
    & - & $\mathcal{C}_0$ & $\mathcal{C}_{0,\mathrm{ott}}$ & $\mathcal{\tau}_0$ & $\mathcal{\tau}_{0,\mathrm{ott}}$ \\
    \hline
%----------------------------------------------------------------------------------------------
    \multicolumn {11}{|c|}{$\alpha=1.2$} \\
    \hline
1000 & 16 &  5 & 5 & 4 &  4   & 11 &  5 & 5 & 4 &  4  \\
2000 & 17 &  5 & 5 & 4 &  4   & 11 &  5 & 5 & 4 &  4  \\
4000 & 17 &  5 & 5 & 4 &  4   & 12 &  5 & 5 & 4 &  4  \\
8000 & 18 &  5 & 5 & 4 &  4   & 12 &  5 & 5 & 4 &  4  \\
   \hline
%----------------------------------------------------------------------------------------------
    \multicolumn {11}{|c|}{$\alpha=1.5$} \\
    \hline
1000 & 57 &  6 & 7 & 4 &  5  & 39 &  6 & 6 & 4 &  4  \\
2000 & 67 &  6 & 7 & 4 &  5  & 45 &  6 & 6 & 4 &  4  \\
4000 & 78 &  6 & 7 & 4 &  5  & 53 &  6 & 6 & 4 &  4  \\
8000 & 92 &  7 & 7 & 4 &  4  & 62 &  6 & 6 & 4 &  4  \\
     \hline
%----------------------------------------------------------------------------------------------
    \multicolumn {11}{|c|}{$\alpha=1.8$} \\
    \hline
1000 & 191 &  7 & 11 & 5 &  5  & 132 &  7 & 9 & 4 &  5  \\
2000 & 249 &  7 & 11 & 5 &  5  & 172 &  7 & 9 & 4 &  5  \\
4000 & 325 &  7 & 11 & 5 &  5  & 225 &  7 & 9 & 5 &  5  \\
8000 & 423 &  8 & 11 & 5 &  5  & 293 &  7 & 9 & 5 &  5  \\
     \hline
  \end{tabular}
  \caption{Number of preconditioned GMRES iterations to solve the linear system with matrix $\mathcal{A}_0^{\mathrm{antiR}}=  \nu I-{k}A_0^{\mathrm{antiR}}$  (Implicit Euler method) and $\mathcal{A}_0^{\mathrm{antiR}}=  \nu I-\frac{k}{2}A_0^{\mathrm{antiR}}$ (Crank-Nicolson method) for increasing dimension $n$  till $tol =1.e-6$ - case $k=1$.} \label{tab:EIcn_gmresit_AntiR}
\end{table}
%%%
\begin{table}
  \centering
  \footnotesize
  \begin{tabular}{|c|ccccc|ccccc|}
  \hline
    &  \multicolumn {5}{c|}{Implicit Euler} &  \multicolumn {5}{c|}{Crank-Nicolson}\\
    \hline
    $n$
    & - & $\mathcal{C}_0$ & $\mathcal{C}_{0,\mathrm{ott}}$ & $\mathcal{\tau}_0$ & $\mathcal{\tau}_{0,\mathrm{ott}}$
    & - & $\mathcal{C}_0$ & $\mathcal{C}_{0,\mathrm{ott}}$ & $\mathcal{\tau}_0$ & $\mathcal{\tau}_{0,\mathrm{ott}}$ \\
    \hline
%----------------------------------------------------------------------------------------------
    \multicolumn {11}{|c|}{$\alpha=1.2$} \\
    \hline
1000 & 197 &  8 & 12 & 7 &  7  & 140 &  8 & 10& 7 &  7  \\
2000 & 214 &  8 & 11 & 5 &  6  & 148 &  8 & 9 & 5 &  5  \\
4000 & 227 &  8 & 10 & 5 &  6  & 156 &  8 & 9 & 5 &  5  \\
8000 & 239 &  9 & 10 & 5 &  5  & 164 &  8 & 9 & 5 &  5  \\
  \hline
%----------------------------------------------------------------------------------------------
    \multicolumn {11}{|c|}{$\alpha=1.5$} \\
    \hline
1000 & 508 &  9 & 21 & 5 &  6  & 401 &  9 & 19 & 5 &  6  \\
2000 & 660 &  9 & 23 & 5 &  6  & 510 &  9 & 19 & 5 &  6  \\
4000 & 832 &  9 & 22 & 5 &  6  & 612 &  9 & 19 & 5 &  6  \\
8000 & 997 &  9 & 22 & 5 &  6  & 725  & 9 & 18 & 5 &  5  \\
     \hline
%----------------------------------------------------------------------------------------------
    \multicolumn {11}{|c|}{$\alpha=1.8$} \\
    \hline
1000 & 846  &  6 & 32 & 5 &  6   & 789  &  6 & 31 & 5 &  6  \\
2000 & 1464 &  6 & 35 & 5 &  6   & 1378 &  6 & 36 & 5 &  6  \\
4000 & 2484 &  6 & 39 & 5 &  6   & 1988 &  6 & 38 & 5 &  6  \\
8000 & 3475 &  6 & 38 & 5 &  6   & 2600 &  6 & 43 & 5 &  6  \\
     \hline
  \end{tabular}
  \caption{Number of preconditioned GMRES iterations to solve the linear system with matrix $\mathcal{A}_0^{\mathrm{antiR}}=  \nu I-{k}A_0^{\mathrm{antiR}}$  (Implicit Euler method) and $\mathcal{A}_0^{\mathrm{antiR}}=  \nu I-\frac{k}{2}A_0^{\mathrm{antiR}}$ (Crank-Nicolson method) for increasing dimension $n$  till $tol =1.e-6$ - case $k=100$.} \label{tab:EIcn_gmresit_AntiR_k100}
\end{table}
%%%%%
\begin{table}
  \centering
  \footnotesize
  \begin{tabular}{|c|ccccc|ccccc|}
  \hline
    &  \multicolumn {5}{c|}{Implicit Euler} &  \multicolumn {5}{c|}{Crank-Nicolson}\\
    \hline
    $n$
    & - & $\mathcal{C}_0$ & $\mathcal{C}_{0,\mathrm{ott}}$ & $\mathcal{\tau}_0$ & $\mathcal{\tau}_{0,\mathrm{ott}}$
    & - & $\mathcal{C}_0$ & $\mathcal{C}_{0,\mathrm{ott}}$ & $\mathcal{\tau}_0$ & $\mathcal{\tau}_{0,\mathrm{ott}}$ \\
    \hline
%----------------------------------------------------------------------------------------------
    \multicolumn {11}{|c|}{$\alpha=1.2$} \\
    \hline
1000 & 21 &  6 & 6 & 4 &  4  & 15 &  5 & 5 & 4 &  4  \\
2000 & 22 &  6 & 6 & 4 &  4  & 16 &  5 & 5 & 4 &  4  \\
4000 & 23 &  6 & 6 & 4 &  4  & 17 &  5 & 5 & 4 &  4  \\
8000 & 24 &  6 & 6 & 4 &  4  & 18 &  5 & 5 & 4 &  4  \\
  \hline
%----------------------------------------------------------------------------------------------
    \multicolumn {11}{|c|}{$\alpha=1.5$} \\
    \hline
1000 & 53 &  7 & 8 & 5 &  5  & 39 &  7 & 7 & 5 &  5  \\
2000 & 59 &  7 & 8 & 5 &  5  & 45 &  7 & 7 & 5 &  5  \\
4000 & 69 &  7 & 8 & 5 &  5  & 52 &  7 & 7 & 5 &  5  \\
8000 & 77 &  8 & 8 & 5 &  5  & 59 &  7 & 7 & 5 &  5  \\
     \hline
%----------------------------------------------------------------------------------------------
    \multicolumn {11}{|c|}{$\alpha=1.8$} \\
    \hline
1000 & 128 &  7 & 13 & 5 &  5 & 98 &  7 & 10 & 5 &  5  \\
2000 & 154 &  7 & 13 & 5 &  5 & 120 &  7 & 10 & 5 &  5  \\
4000 & 192 &  8 & 12 & 5 &  5 & 149 &  7 & 10 & 5 &  5  \\
8000 & 217 &  6 & 12 & 5 &  5 & 172 &  7 & 10 & 5 &  5  \\
     \hline
  \end{tabular}
  \caption{Number of preconditioned GMRES iterations to solve the linear system with matrix $\mathcal{A}_0^{\mathrm{antiR}}=  \nu I-{k}A_0^{\mathrm{antiR}}$  (Implicit Euler method) and $\mathcal{A}_0^{\mathrm{antiR}}=  \nu I-\frac{k}{2}A_0^{\mathrm{antiR}}$ (Crank-Nicolson method) for increasing dimension $n$  till $tol =1.e-6$ - case $k=1$  random exact solution.} \label{tab:EIcn_gmresit_AntiR_random}
\end{table}
%%
%----------------------------------------------------------------------------------------------
\section{Truncated approximations, the anti-reflective transform, and numerical experiments}\label{sec:truncations}

In order to increase the computational efficiency we may consider a differently truncated version of the previous approximation of the left and right
fractional Riemann-Liouville derivatives, suitable tailored so that the arising matrix belongs to the anti-reflective matrix algebra \cite{paperAR} and the solution of the linear systems can be achieved by a direct solver within $O(N\log N)$ real arithmetic operations via few fast discrete sine transform of type I
(for other fast trigonometric and Fourier-like transforms and their use, see \cite{DiB-S,KO} and references there reported).
As emphasized in \cite{Van} the cost of one fast discrete transform of type I is around one half of the cost of the celebrated fast Fourier transform.
Furthermore the related solver is of direct type and we do not need any preconditioned Krylov iterative solver so that the overall cost is much lower, when compared with the techniques proposed for the nontruncated versions.

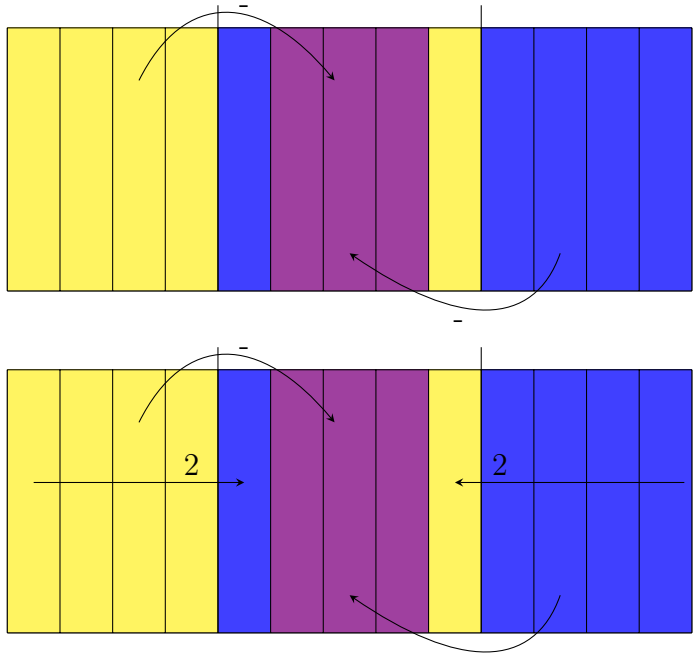
\begin{figure}
%Anti-symmetric BCs
\begin{center}
\begin{tikzpicture}{0cm}{0cm}{7.5cm}{5cm}
%-----------------------------------------------------------------------
% Griglia
\pgfgrid[stepx=0.7cm,stepy=3.5cm]{\pgfxy(0,0)}{\pgfxy(9.1,3.5)} %griglia totale
%%%%%%%%%%%%%%%%%%%%%%%%%%%%%%%%%%%%%%%%%%%%%%%%%%%%%%%%%%%%%%%%%%%%%%%%%%%%%%%%%%%%%%           %-----------------------------------------------------------------------
% Rettangolo esterno sinistra
\begin{colormixin}{75!white} \color{yellow}
\pgfrect[fill]{\pgfxy(0,0)}{\pgfxy(2.8,3.5)}
\pgffill
\end{colormixin}
%-----------------------------------------------------------------------
% Rettangolo interno sinistra
\begin{colormixin}{75!white} \color{yellow}
\pgfrect[fill]{\pgfxy(3.5,0)}{\pgfxy(2.8,3.5)}
\pgffill
\end{colormixin}
%-----------------------------------------------------------------------
% Rettangolo esterno destra
\begin{colormixin}{75!white} \color{blue}
\pgfrect[fill]{\pgfxy(6.3,0)}{\pgfxy(2.8,3.5)}
\pgffill
\end{colormixin}
%-----------------------------------------------------------------------
% Rettangolo basso interno
\begin{colormixin}{75!white} \color{blue}
\pgfrect[fill]{\pgfxy(2.8,0)}{\pgfxy(2.8,3.5)}
\pgffill
\end{colormixin}
%-----------------------------------------------------------------------
% Rettangolo centrale interno
\begin{colormixin}{75!white} \color{violet}
\pgfrect[fill]{\pgfxy(3.5,0)}{\pgfxy(2.1,3.5)}
\pgffill
\end{colormixin}
%*********************************************************************************
% Linee verticali
\color{black}
\pgfxyline(2.8,0)(2.8,3.8)
\pgfxyline(6.3,0)(6.3,3.8)
%*********************************************************************************
% Linee diagonali
%\pgfxyline(2.8,3.5)(6.3,0) %linee diagonali psf
%\pgfxyline(0.7,3.5)(4.2,0) %linee diagonali psf
%\pgfxyline(4.9,3.5)(8.4,0) %linee diagonali psf
%*********************************************************************************
\pgfsetarrowsend{stealth}
\pgfxycurve(1.75,2.8)(2.35,4)(3.35,4)(4.35,2.8) %curva centro immagine/centro psf
\pgfputat{\pgfxy(3.15,3.7)}{\pgfbox[center,base]{-}} %scritta psf
\pgfclearendarrow
\pgfsetstartarrow{\pgfarrowtriangle{4pt}} % stile punta freccia per curva
\pgfxycurve(7.35,0.5)(7,-0.5)(6,-0.5)(4.55,0.5) %curva centro immagine/centro psf
\pgfputat{\pgfxy(6,-0.5)}{\pgfbox[center,base]{-}} %scritta psf
\pgfclearstartarrow
\pgfclearendarrow
%*********************************************************************************
% Ridisegna contorno
\pgfgrid[stepx=0.7cm,stepy=3.5cm]{\pgfxy(0,0)}{\pgfxy(9.1,3.5)} %griglia totale
%*********************************************************************************
\end{tikzpicture}
%\end{center}
% ----------------------------------------------------------------
% ----------------------------------------------------------------
%Anti-Reflective BCs
%\begin{center}
\begin{tikzpicture}{0cm}{0cm}{7.5cm}{5cm}
%-----------------------------------------------------------------------
% Griglia
\pgfgrid[stepx=0.7cm,stepy=3.5cm]{\pgfxy(0,0)}{\pgfxy(9.1,3.5)} %griglia totale
%%%%%%%%%%%%%%%%%%%%%%%%%%%%%%%%%%%%%%%%%%%%%%%%%%%%%%%%%%%%%%%%%%%%%%%%%%%%%%%%%%%%%%           %-----------------------------------------------------------------------
% Rettangolo esterno sinistra
\begin{colormixin}{75!white} \color{yellow}
\pgfrect[fill]{\pgfxy(0,0)}{\pgfxy(2.8,3.5)}
\pgffill
\end{colormixin}
%-----------------------------------------------------------------------
% Rettangolo interno sinistra
\begin{colormixin}{75!white} \color{yellow}
\pgfrect[fill]{\pgfxy(3.5,0)}{\pgfxy(2.8,3.5)}
\pgffill
\end{colormixin}
%-----------------------------------------------------------------------
% Rettangolo esterno destra
\begin{colormixin}{75!white} \color{blue}
\pgfrect[fill]{\pgfxy(6.3,0)}{\pgfxy(2.8,3.5)}
\pgffill
\end{colormixin}
%-----------------------------------------------------------------------
% Rettangolo basso interno
\begin{colormixin}{75!white} \color{blue}
\pgfrect[fill]{\pgfxy(2.8,0)}{\pgfxy(2.8,3.5)}
\pgffill
\end{colormixin}
%-----------------------------------------------------------------------
% Rettangolo centrale interno
\begin{colormixin}{75!white} \color{violet}
\pgfrect[fill]{\pgfxy(3.5,0)}{\pgfxy(2.1,3.5)}
\pgffill
\end{colormixin}
%*********************************************************************************
% Linee verticali
\color{black}
\pgfxyline(2.8,0)(2.8,3.8)
\pgfxyline(6.3,0)(6.3,3.8)
%*********************************************************************************
% Linee diagonali
%\pgfxyline(2.8,3.5)(6.3,0) %linee diagonali psf
%\pgfxyline(0.7,3.5)(4.2,0) %linee diagonali psf
%\pgfxyline(4.9,3.5)(8.4,0) %linee diagonali psf
%*********************************************************************************
\pgfsetarrowsend{stealth}
\pgfxycurve(1.75,2.8)(2.35,4)(3.35,4)(4.35,2.8) %curva centro immagine/centro psf
\pgfputat{\pgfxy(3.15,3.7)}{\pgfbox[center,base]{-}} %scritta psf
\pgfclearendarrow
\pgfsetstartarrow{\pgfarrowtriangle{4pt}} % stile punta freccia per curva
\pgfxycurve(7.35,0.5)(7,-0.5)(6,-0.5)(4.55,0.5) %curva centro immagine/centro psf
\pgfputat{\pgfxy(6,-0.5)}{\pgfbox[center,base]{-}} %scritta psf
\pgfclearstartarrow
\pgfclearendarrow
%*********************************************************************************
\pgfxyline(0.35,2)(3.15,2) %curva centro immagine/centro psf
\pgfputat{\pgfxy(2.45,2.1)}{\pgfbox[center,base]{2}} %scritta psf
\pgfxyline(9,2)(5.95,2) %curva centro immagine/centro psf
\pgfputat{\pgfxy(6.55,2.1)}{\pgfbox[center,base]{2}} %scritta psf
\pgfclearendarrow
%*********************************************************************************
% Ridisegna contorno
\pgfgrid[stepx=0.7cm,stepy=3.5cm]{\pgfxy(0,0)}{\pgfxy(9.1,3.5)} %griglia totale
%*********************************************************************************
\end{tikzpicture}
\end{center}
\vskip -1cm
\caption{Effect of the numerical anti-symmetric and anti-reflective BCs on the matrix structure.} \label{fig:scheme}
\end{figure}

\begin{figure}\centering %Truncated Anti-symmetric BCs
\begin{center}
\begin{tikzpicture}{0cm}{0cm}{7.5cm}{5cm}
%-----------------------------------------------------------------------
% Griglia
\pgfgrid[stepx=0.7cm,stepy=3.5cm]{\pgfxy(0,0)}{\pgfxy(9.1,3.5)} %griglia totale
%%%%%%%%%%%%%%%%%%%%%%%%%%%%%%%%%%%%%%%%%%%%%%%%%%%%%%%%%%%%%%%%%%%%%%%%%%%%%%%%%%%%%%           %-----------------------------------------------------------------------
% Triangolo esterno sinistra
\begin{colormixin}{75!white} \color{yellow}
\pgfmoveto{\pgfxy(0.7,3.5)}
\pgflineto{\pgfxy(2.8,3.5)}
\pgflineto{\pgfxy(2.8,1.4)}
\pgflineto{\pgfxy(0.7,3.5)}
\pgffill
\end{colormixin}
%-----------------------------------------------------------------------
% Triangolo interno sinistra
\begin{colormixin}{75!white} \color{yellow}
\pgfmoveto{\pgfxy(3.5,3.5)}
\pgflineto{\pgfxy(5.6,3.5)}
\pgflineto{\pgfxy(3.5,1.4)}
\pgflineto{\pgfxy(3.5,3.5)}
\pgffill
\end{colormixin}
%-----------------------------------------------------------------------
% Triangolo esterno destra
\begin{colormixin}{75!white} \color{blue}
\pgfmoveto{\pgfxy(6.3,2.1)}
\pgflineto{\pgfxy(8.4,0)}
\pgflineto{\pgfxy(6.3,0)}
\pgflineto{\pgfxy(6.3,2.1)}
\pgffill
\end{colormixin}
%-----------------------------------------------------------------------
% Triangolo basso interno
\begin{colormixin}{75!white} \color{blue}
\pgfmoveto{\pgfxy(3.5,0)}
\pgflineto{\pgfxy(5.6,2.1)}
\pgflineto{\pgfxy(5.6,0)}
\pgflineto{\pgfxy(3.5,0)}
\pgffill
\end{colormixin}
%*********************************************************************************
% Linee verticali
\color{black}
\pgfxyline(2.8,0)(2.8,3.8)
\pgfxyline(6.3,0)(6.3,3.8)
%*********************************************************************************
% Linee diagonali
\pgfxyline(2.8,3.5)(6.3,0) %linee diagonali psf
\pgfxyline(0.7,3.5)(4.2,0) %linee diagonali psf
\pgfxyline(4.9,3.5)(8.4,0) %linee diagonali psf
%*********************************************************************************
\pgfsetarrowsend{stealth}
\pgfxycurve(1.75,2.8)(2.35,4)(3.35,4)(4.35,2.8) %curva centro immagine/centro psf
\pgfputat{\pgfxy(3.15,3.7)}{\pgfbox[center,base]{-}} %scritta psf
\pgfclearendarrow
\pgfsetstartarrow{\pgfarrowtriangle{4pt}} % stile punta freccia per curva
\pgfxycurve(7.35,0.5)(7,-0.5)(6,-0.5)(4.55,0.5) %curva centro immagine/centro psf
\pgfputat{\pgfxy(6,-0.5)}{\pgfbox[center,base]{-}} %scritta psf
\pgfclearstartarrow
\pgfclearendarrow
%*********************************************************************************
% Ridisegna contorno
\pgfgrid[stepx=0.7cm,stepy=3.5cm]{\pgfxy(0,0)}{\pgfxy(9.1,3.5)} %griglia totale
%*********************************************************************************
\end{tikzpicture}
%\end{center}
% ----------------------------------------------------------------
% ----------------------------------------------------------------
%\ \\ Truncated Anti-Reflective BCs
%\begin{center}
\begin{tikzpicture}{0cm}{0cm}{7.5cm}{5cm}
%-----------------------------------------------------------------------
% Griglia
\pgfgrid[stepx=0.7cm,stepy=3.5cm]{\pgfxy(0,0)}{\pgfxy(9.1,3.5)} %griglia totale
%%%%%%%%%%%%%%%%%%%%%%%%%%%%%%%%%%%%%%%%%%%%%%%%%%%%%%%%%%%%%%%%%%%%%%%%%%%%%%%%%%%%%%           %-----------------------------------------------------------------------
% Triangolo esterno sinistra
\begin{colormixin}{75!white} \color{yellow}
\pgfmoveto{\pgfxy(0.7,3.5)}
\pgflineto{\pgfxy(2.8,3.5)}
\pgflineto{\pgfxy(2.8,1.4)}
\pgflineto{\pgfxy(0.7,3.5)}
\pgffill
\end{colormixin}
%-----------------------------------------------------------------------
% Triangolo interno sinistra
\begin{colormixin}{75!white} \color{yellow}
\pgfmoveto{\pgfxy(3.5,3.5)}
\pgflineto{\pgfxy(5.6,3.5)}
\pgflineto{\pgfxy(3.5,1.4)}
\pgflineto{\pgfxy(3.5,3.5)}
\pgffill
\end{colormixin}
%-----------------------------------------------------------------------
% Triangolo esterno destra
\begin{colormixin}{75!white} \color{blue}
\pgfmoveto{\pgfxy(6.3,2.1)}
\pgflineto{\pgfxy(8.4,0)}
\pgflineto{\pgfxy(6.3,0)}
\pgflineto{\pgfxy(6.3,2.1)}
\pgffill
\end{colormixin}
%-----------------------------------------------------------------------
% Triangolo basso interno
\begin{colormixin}{75!white} \color{blue}
\pgfmoveto{\pgfxy(3.5,0)}
\pgflineto{\pgfxy(5.6,2.1)}
\pgflineto{\pgfxy(5.6,0)}
\pgflineto{\pgfxy(3.5,0)}
\pgffill
\end{colormixin}
%-----------------------------------------------------------------------
% Rettangolo interno sinistra
\begin{colormixin}{75!white} \color{gray}
\pgfrect[fill]{\pgfxy(2.8,0)}{\pgfxy(0.7,3.5)}
\pgffill
\end{colormixin}
%-----------------------------------------------------------------------
% Rettangolo interno destro
\begin{colormixin}{75!white} \color{gray}
\pgfrect[fill]{\pgfxy(5.6,0)}{\pgfxy(0.7,3.5)}
\pgffill
\end{colormixin}
%*********************************************************************************
% Linee verticali
\color{black}
\pgfxyline(2.8,0)(2.8,3.8)
\pgfxyline(6.3,0)(6.3,3.8)
%*********************************************************************************
% Linee diagonali
\pgfxyline(2.8,3.5)(6.3,0) %linee diagonali psf
\pgfxyline(0.7,3.5)(4.2,0) %linee diagonali psf
\pgfxyline(4.9,3.5)(8.4,0) %linee diagonali psf
%*********************************************************************************
\pgfsetarrowsend{stealth}
\pgfxycurve(1.75,2.8)(2.35,4)(3.35,4)(4.35,2.8) %curva centro immagine/centro psf
\pgfputat{\pgfxy(3.15,3.7)}{\pgfbox[center,base]{-}} %scritta psf
\pgfclearendarrow
\pgfsetstartarrow{\pgfarrowtriangle{4pt}} % stile punta freccia per curva
\pgfxycurve(7.35,0.5)(7,-0.5)(6,-0.5)(4.55,0.5) %curva centro immagine/centro psf
\pgfputat{\pgfxy(6,-0.5)}{\pgfbox[center,base]{-}} %scritta psf
\pgfclearstartarrow
%*********************************************************************************
\pgfxyline(1.75,2.8)(3.15,2.8) %curva centro immagine/centro psf
\pgfputat{\pgfxy(2.45,2.4)}{\pgfbox[center,base]{2}} %scritta psf
\pgfxyline(7.35,0.5)(5.95,0.5) %curva centro immagine/centro psf
\pgfputat{\pgfxy(6.55,0.6)}{\pgfbox[center,base]{2}} %scritta psf
\pgfclearendarrow
%*********************************************************************************
% Ridisegna contorno
\pgfgrid[stepx=0.7cm,stepy=3.5cm]{\pgfxy(0,0)}{\pgfxy(9.1,3.5)} %griglia totale
%*********************************************************************************
\end{tikzpicture}
\end{center}
\vskip -1cm
\caption{Truncated anti-symmetric and anti-reflective boundaries effect on the matrix structure.} \label{fig:truncated_scheme}
\end{figure}
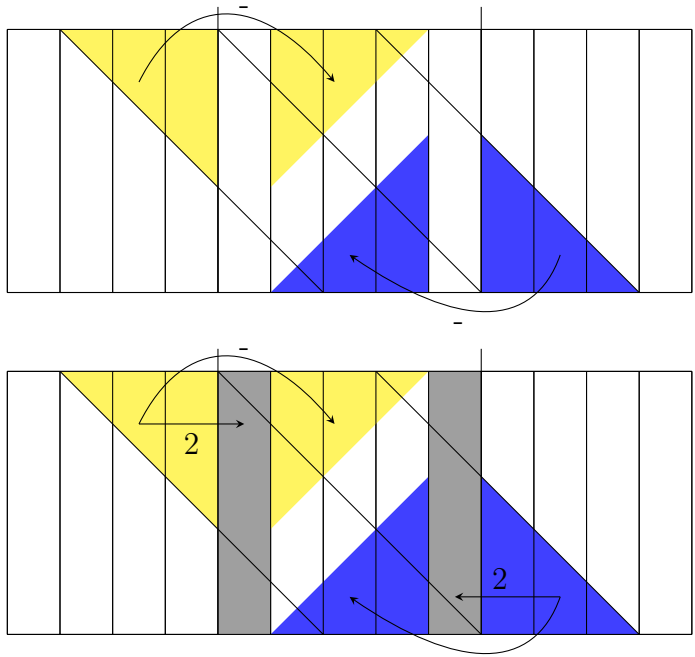

More precisely, we will consider the approximations as follows
\begin{eqnarray*}
 {}_{-\infty}^{RL}D_{x}^{\alpha,ref} u (x_j,t_n)
 &\approx
 & \frac{1}{(\Delta x)^\alpha}\sum_{k=0}^{j+1} g_{k}^\alpha U_{j+1-k}^n
 +\frac{1}{(\Delta x)^\alpha}\sum_{k=j+2}^{N} g_{k}^\alpha U_{j+1-k}^n,\\
 {}_{x}^{RL}D_{\infty}^{\alpha,ref} u (x_j,t_n)
 &\approx & \frac{1}{(\Delta x)^\alpha}\sum_{k=0}^{N-j+1} g_{k}^\alpha U_{j-1+k}^n
 +\frac{1}{(\Delta x)^\alpha}\sum_{k=N-j+2}^{N} g_{k}^\alpha U_{j-1+k}^n,
\end{eqnarray*}
where the number of terms approximating the derivatives in each point is constant.
Thus, before imposing the considered numerical BCs, the resulting matrix structure of the linear system is as follows
\[
A_0^{\mathrm{anti,full}} {\bf U}^{n,\mathrm{full}} = {\bf b}
\]
with ${\bf U}^{n,\mathrm{full}}=[ U_{-N}^n, \dots, U_{-1}^n |U_{0}^n, \dots, U_{N}^n| U_{N+1}^n, \dots, U_{2N}^n]^T$ and
$2A_0^{\mathrm{anti,full}}$ taking the form below
\[\scriptsize
\phantom{+} \left [\begin{array}{ccc ccc | cccc ccc | ccc ccc}
%---------------------------------------------------------           %
0 & \g{N}   & \g{N-1}& \ldots  & \g{3}  & \g{2}  &\g{1}   & \g{0}  & 0       &        &         &        &       &\phantom{\g{2}}   & \phantom{\g{3}}  & \phantom{\ldots} & \phantom{\g{N-1}}  & \phantom{\g{N}} & \phantom{0}   \\
  &  0      & \g{N}  & \g{N-1} & \ldots & \g{3}  &\g{2}   & \g{1}  & \g{0}   & 0      &         &        &       &      &   & & &  &   \\
  &         & \ddots & \ddots  & \ddots & \ddots &\ddots  & \ddots & \ddots  & \ddots & \ddots  &        &       &      &   & & &  &   \\
  &         &        & \ddots  & \ddots & \ddots &\ddots  & \ddots & \ddots  & \ddots & \ddots  & \ddots &       &      &   & & &  &   \\
  &         &        &         & \ddots & \ddots &\ddots  & \ddots & \ddots  & \ddots & \ddots  & \ddots &       &      &   & & &  &   \\
  &         &        &         &        &   0    &\g{N}   & \g{N-1}& \ldots  & \ldots & \ldots  & \g{1}  & \g{0} & 0    &   & & &  &   \\
  &         &        &         &        &        &   0    & \g{N}  & \g{N-1} & \ldots & \g{3}   & \g{2}  & \g{1} &\g{0} & 0 & & &  &   \\
\end{array}\right] \\
\]
\[\scriptsize
+ \quad\left [\begin{array}{ccc ccc | cccc ccc | ccc ccc}
%---------------------------------------------------------
& & & & 0 & \g{0} & \g{1}   & \g{2}  & \g{3}   & \ldots & \g{N-1} & \g{N}   & 0      &        &        &        &          &       &\\
& & & &   & 0     & \g{0}   & \g{1}  & \g{2}   & \g{3}  &  \ldots & \g{N-1} & \g{N}  &0       &        &        &          &       & \\
& & & &   &       & \ddots  & \ddots & \ddots  & \ddots & \ddots  & \ddots  & \ddots &\ddots  & \ddots &        &          &       & \\
& & & &   &       &         & \ddots & \ddots  & \ddots & \ddots  & \ddots  & \ddots &\ddots  & \ddots & \ddots &          &       &\\
& & & &   &       &         &        & \ddots  & \ddots & \ddots  & \ddots  & \ddots &\ddots  & \ddots & \ddots & \ddots   &       & \\
& & & &   &       &         &        &         &  0     & \g{0}   & \g{1}   & \g{2}  &\g{3}   & \ldots & \g{N-1}& \g{N}    & 0     & \\
\phantom{0} & \phantom{\g{N}}   & \phantom{\g{N-1}}& \phantom{\ldots}  & \phantom{\g{3}}  & \phantom{\g{2}}&         &        &         &        & 0       & \g{0}   & \g{1}  &\g{2}   & \g{3}  & \ldots & \g{N-1}  & \g{N} & 0\\
\end{array}\right].
\]
Thus, by imposing the anti-symmetric BCs, we square the system and the matrix becomes
\[ \scriptsize
\frac{1}{2}\begin{bmatrix}
\begin{array}{c|cccccc|c}
2\g{1}      & \g{0}+\g{2}-\textcolor{brown}{(\g{0}+\g{2})}  & \g{3}-\textcolor{brown}{\g{3}} &  \ldots &  \ldots &  \g{N-1}-\textcolor{brown}{\g{N-1}}  & \g{N}-\textcolor{brown}{\g{N}}  & 0\\
\hline
\g{0}+\g{2} & 2\g{1}-\textcolor{brown}{\g{3}}      & \g{0}+\g{2}-\textcolor{brown}{\g{4}} & \g{3}-\textcolor{brown}{\g{5}}  &        & \ldots & \g{N-1}                         &  \g{N}      \\
\g{3}       & \g{0}+\g{2}-\textcolor{brown}{\g{4}} & \vdots      & \vdots & \vdots & \vdots & \g{N-2}-\textcolor{brown}{\g{N}}                         & \g{N-1}     \\
\vdots      & \vdots                               & \vdots      & \vdots & \vdots & \vdots & \vdots                          & \vdots      \\
\g{N-1}     & \g{N-2}-\textcolor{brown}{\g{N}}     & \vdots                          & \vdots & \vdots & \vdots & \vdots    & \g{3}       \\
\g{N}       & \g{N-1}                              & \g{N-2}-\textcolor{brown}{\g{N}}     & \ldots & \ldots & \ldots & 2\g{1} -\textcolor{brown}{\g{3}}& \g{0}+\g{2} \\
\hline
   0        & \g{N}-\textcolor{brown}{\g{N}}              & \g{N-1}-\textcolor{brown}{\g{N-1}}     & \ldots & \ldots  & \g{3}    & \g{0}+\g{2} -\textcolor{brown}{(\g{0}+\g{2})}  & 2\g{1}
\end{array}
\end{bmatrix},
\]
that is the anti-symmetric matrix
\begin{equation}\label{structure-anti-symm}
\frac{1}{2}\begin{bmatrix}
\begin{array}{c|cccccc|c}
2\g{1}      & 0  & 0 &  \ldots &  \ldots &  0  & 0  & 0\\
\hline
\g{0}+\g{2} &             &       &  &  &  &  & \g{N}      \\
\g{3}       &             &       &  &  &  &  & \g{N-1}     \\
\vdots      &             &       & \tau(t_0,t_1,t_2, \ldots,t_N) &  &  &  & \vdots      \\
\g{N-1}     &             &       &  &  &  &  & \g{3}       \\
\g{N}       &             &       &  &  &  &  & \g{0}+\g{2} \\
\hline
   0        & 0              & 0     & \ldots & \ldots  & 0    & 0  & 2\g{1}
\end{array}
\end{bmatrix},
\end{equation}
where $\tau(t_0,t_1,t_2, \ldots,t_N)$ is the matrix belonging to the $\tau$ algebra (associated to the discrete sine transform of type I) with coefficients $\{t_i\}_{i=0,\ldots,N}$ given by the Toeplitz coefficients defined in \eqref{eq:Tcoeff}.
In the same way the matrix obtained by imposing numerical anti-reflective BCs shows the expression
\begin{equation}\label{structure-anti-refl}
\frac{1}{2}\begin{bmatrix}
\begin{array}{c|cccccc|c}
2\g{1} \textcolor{blue}{+ 2\g{0}} + \tilde{z}_1     & 0  & 0 &  \ldots &  \ldots &  0  & 0  & 0\\
\hline
\g{0}+\g{2}+ \tilde{z}_2  &       &       &  &  &  &  & \g{N}      \\%+ \tilde{z}_{N}
\g{3} + \tilde{z}_3       &       &       &  &  &  &  & \g{N-1} + \tilde{z}_{N-1}     \\
\vdots      &             &       & \tau(t_0,t_1,t_2, \ldots,t_N) &  &  &  & \vdots      \\
\g{N-1}+ \tilde{z}_{N-1}  &       &       &  &  &  &  & \g{3} + \tilde{z}_{3}      \\
\g{N}%+ \tilde{z}_{N}
&       &       &  &  &  &  & \g{0}+\g{2} + \tilde{z}_{2}\\
\hline
   0        & 0              & 0     & \ldots & \ldots  & 0    & 0  & 2\g{1} \textcolor{blue}{+ 2\g{0}} + \tilde{z}_1
\end{array}
\end{bmatrix},
\end{equation}
where
\[
\tilde{z}_{r} =2 \sum_{k=r+1}^{N} \g{k}, \quad r=1, \ldots, N-1
\]
are just the truncated version of the previously defined quantities $z_{r}$.\\
Therefore, the choice between the two different types of conditions pertains to the quality of the approximation and does not influence the computational efficiency. \\
In Figures \ref{fig:scheme} and \ref{fig:truncated_scheme} the differences in applying the considered numerical BCs in full and its truncated version on $A_0^{\mathrm{anti,full}}$ is highlighted. Notice as the external yellow and blue full wings in Figure \ref{fig:scheme} are transferred inside the matrix $A_0^{\mathrm{anti}}$ with a purple overlap in the central part in the first case, while no overlap is originated in the truncated one.\\
Furthermore, when using physical Dirichlet BCs, the first and the last rows and columns in (\ref{structure-anti-symm}), (\ref{structure-anti-refl}) just disappear and the resulting linear systems become even simpler, since the coefficient matrix reduces to $\tau(t_0,t_1,t_2, \ldots,t_N)$.
Therefore the whole computational cost is that of three sine I transforms owing to the $\tau$ structure \cite{BC} of $\tau(t_0,t_1,t_2, \ldots,t_N)$.
%----------------------------------------------------------------------------------------------
% ----------------------------------------------------------------
%&&&&&&&&&&&&&&&&&&&&&&&&&&&
% ----------------------------------------------------------------
% ----------------------------------------------------------------
% ----------------------------------------------------------------
\subsection*{Numerical precision}
%\label{sec:numprec}
%------------------------------------------------
Finally, we want to check the numerical precision of numerical anti-symmetric and anti-reflective BCs, by considering both nontruncated and truncated versions.
Comparison is made also with respect to more classical numerical BCs of  Dirichlet and reflective (Neumann) type. Notice that the numerical reflective BCs also leads to an algebra of matrices related to fast transform and hence the computational cost is essentially the same as that arising with numerical anti-reflective BCs: however the numerical anti-reflective BCs are more accurate as discussed in \cite{S-AR-proposal,imaging-reflective} and references therein.
\par
\begin{figure}[h]
{\center
\includegraphics[width=\textwidth]{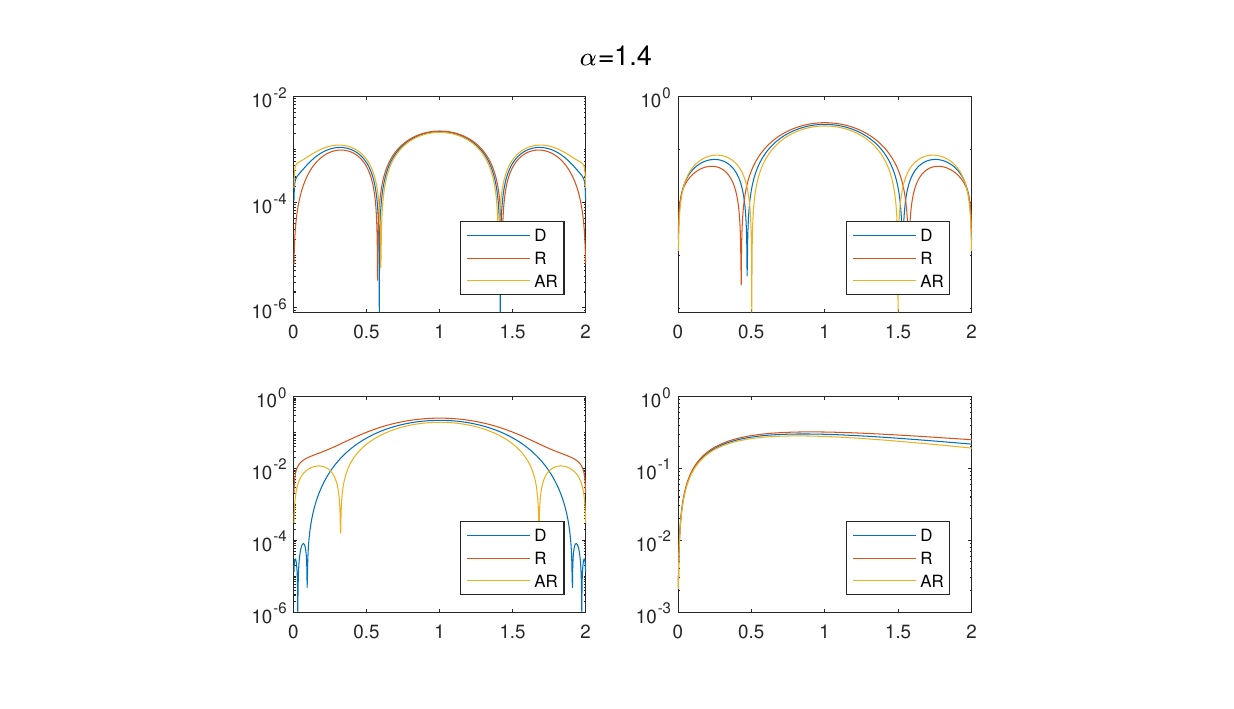}
}
\vskip -1cm \caption{Absolute error vs $x$ at $t \in \{2.e-3,1,2\}$ and Maximal absolute error vs $t$ - Implicit Euler,  $\alpha=1.4$ - Homogeneus Dirichlet BCs.}
\label{fig:EI_14x5_omo}
\end{figure}
\begin{figure}[h]
{\center
\includegraphics[width=\textwidth]{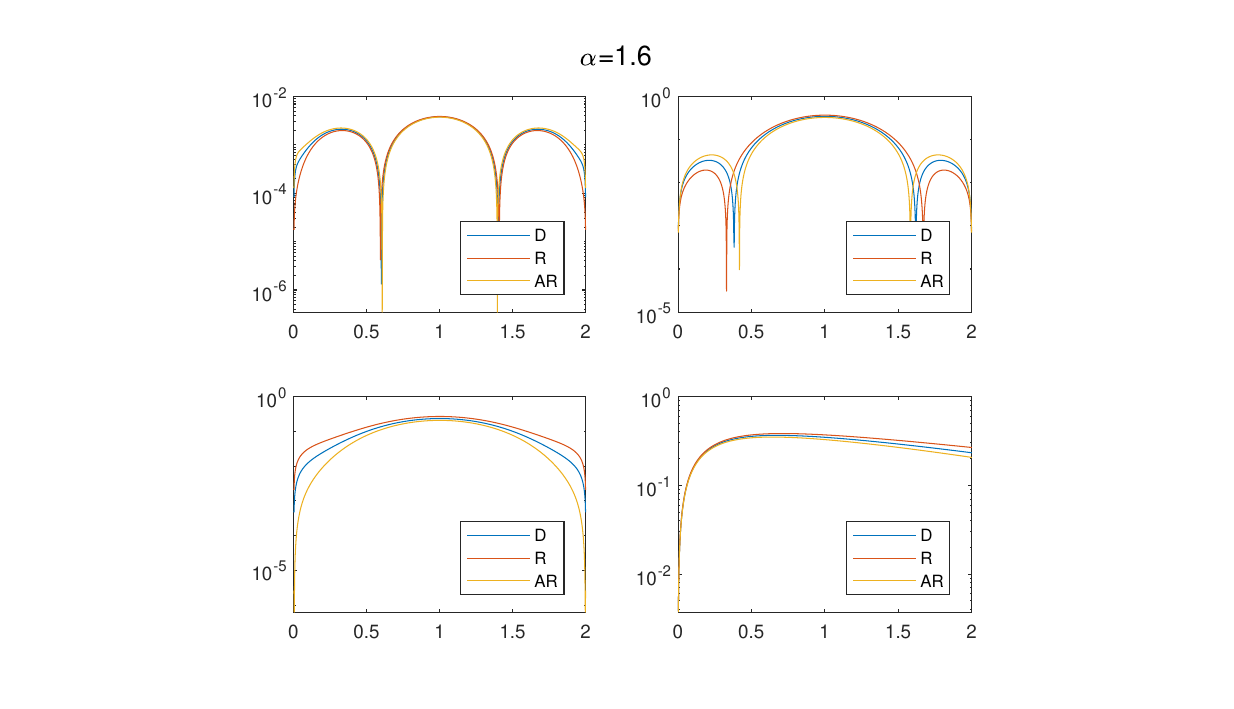}
}
\vskip -1cm \caption{Absolute error vs $x$ at $t \in \{2.e-3,1,2\}$ and Maximal absolute error vs $t$ - Implicit Euler,  $\alpha=1.6$ - Homogeneus Dirichlet BCs.}
\label{fig:EI_16x5_omo}
\end{figure}
\begin{figure}[h]
{\center
\includegraphics[width=\textwidth]{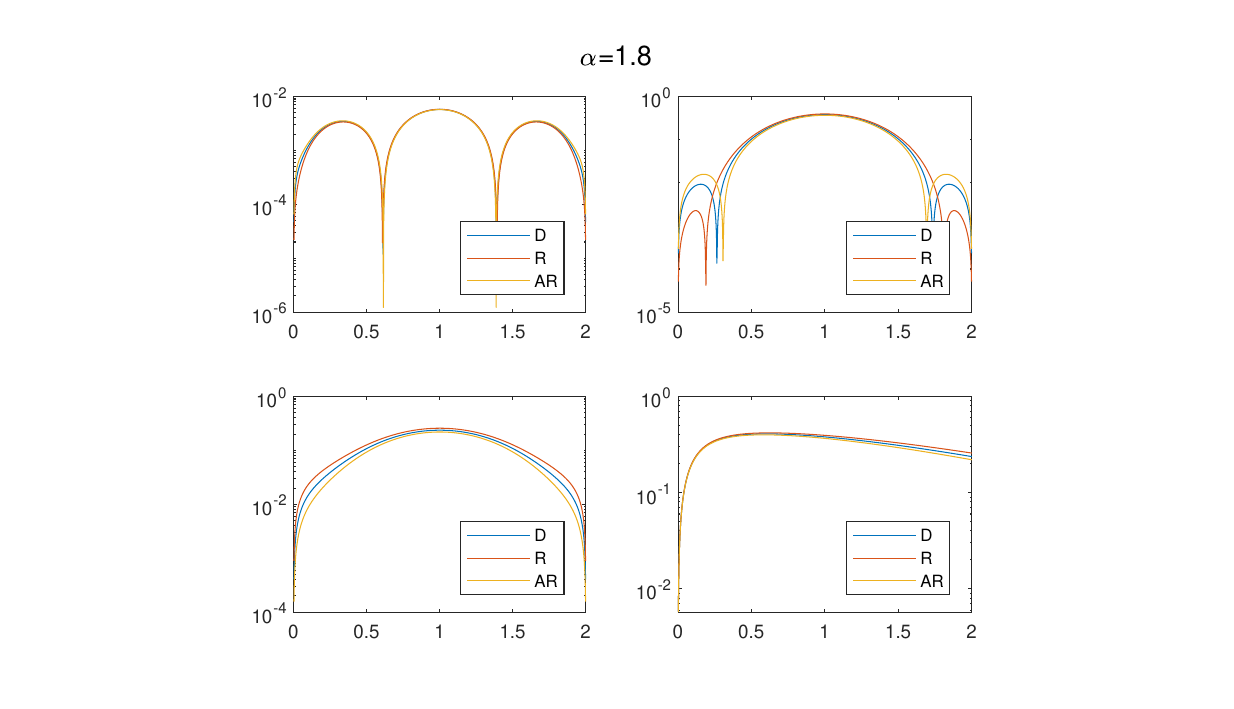}
}
\vskip -1cm \caption{Absolute error vs $x$ at $t \in \{2.e-3,1,2\}$ and Maximal absolute error vs $t$ - Implicit Euler,  $\alpha=1.8$ - Homogeneus Dirichlet BCs.}
\label{fig:EI_18x5_omo}
\end{figure}
\begin{figure}[h]
{\center
\includegraphics[width=\textwidth]{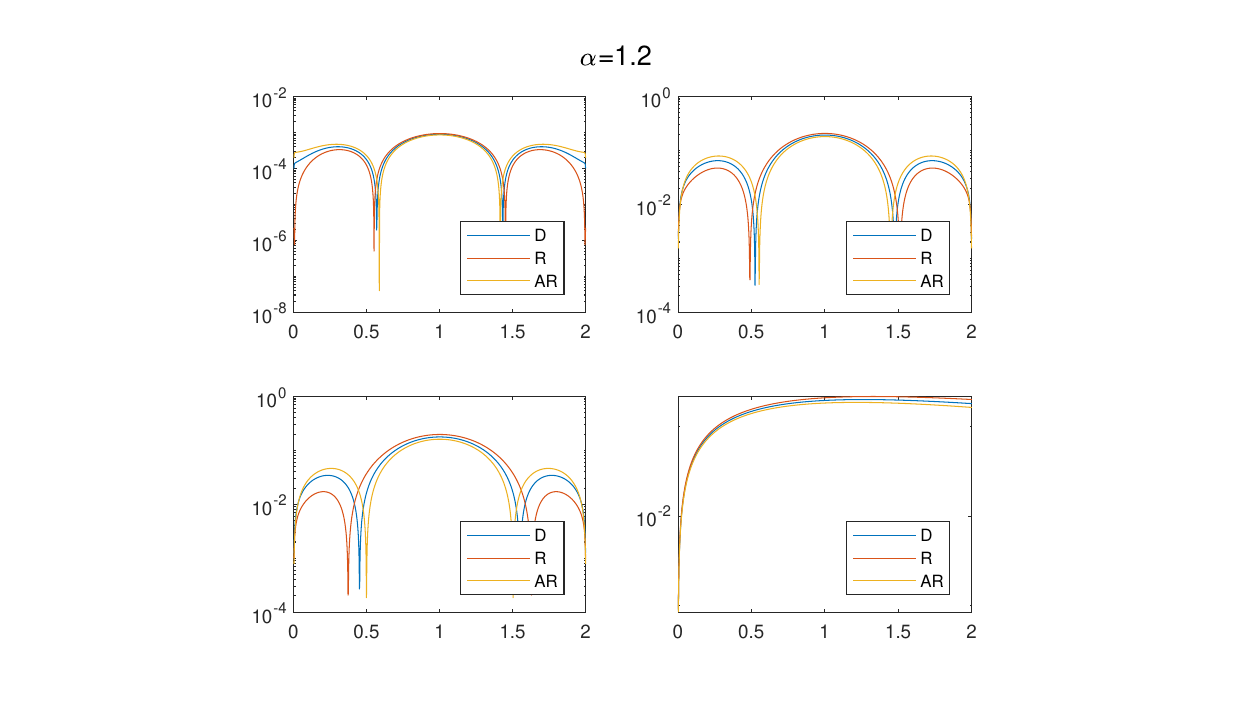}
}
\vskip -1cm \caption{Absolute error vs $x$ at $t \in \{2.e-3,1,2\}$ and Maximal absolute error vs $t$ - Crank-Nicolson,  $\alpha=1.2$ - Homogeneus Dirichlet BCs.}
\label{fig:CN_12x5_omo}
\end{figure}
\begin{figure}[h]
{\center
\includegraphics[width=\textwidth]{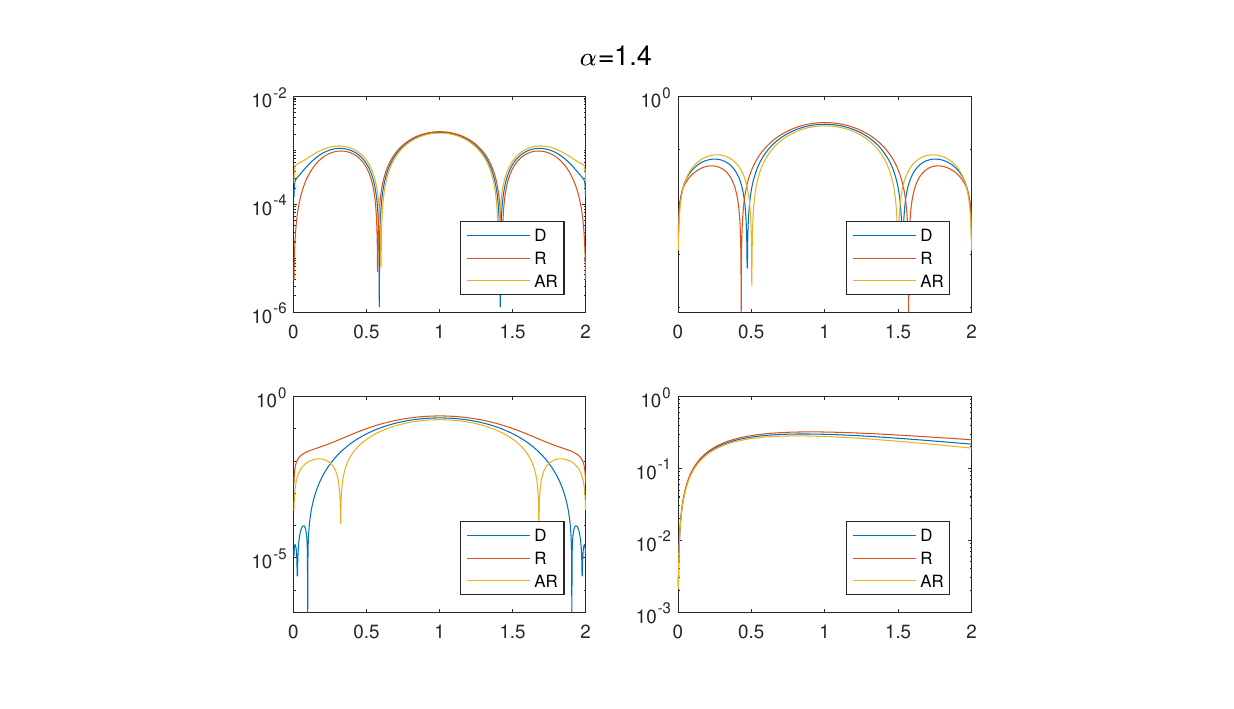}
}
\vskip -1cm \caption{Absolute error vs $x$ at $t \in \{2.e-3,1,2\}$ and Maximal absolute error vs $t$ - Crank-Nicolson,  $\alpha=1.4$ - Homogeneus Dirichlet BCs.}
\label{fig:CN_14x5_omo}
\end{figure}
\begin{figure}[h]
{\center
\includegraphics[width=\textwidth]{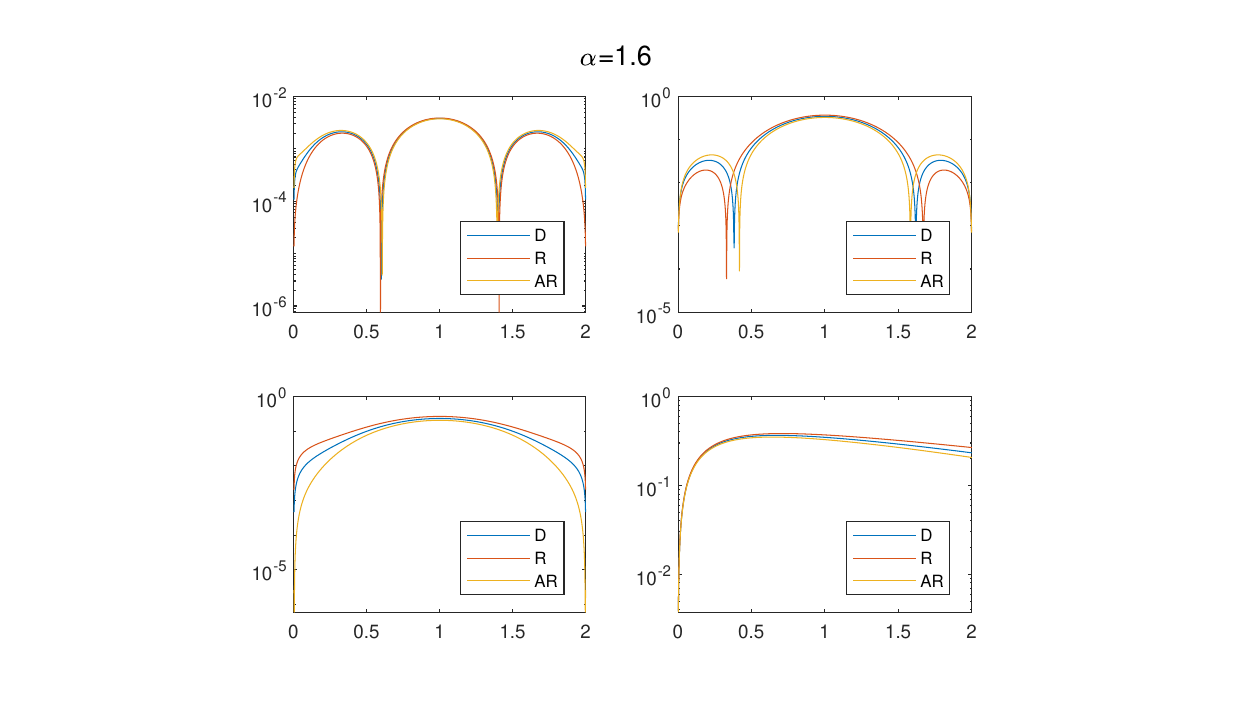}
}
\vskip -1cm \caption{Absolute error vs $x$ at $t \in \{2.e-3,1,2\}$ and Maximal absolute error vs $t$ - Crank-Nicolson,  $\alpha=1.6$ - Homogeneus Dirichlet BCs.}
\label{fig:CN_16x5_omo}
\end{figure}
\begin{figure}[h]
{\center
\includegraphics[width=\textwidth]{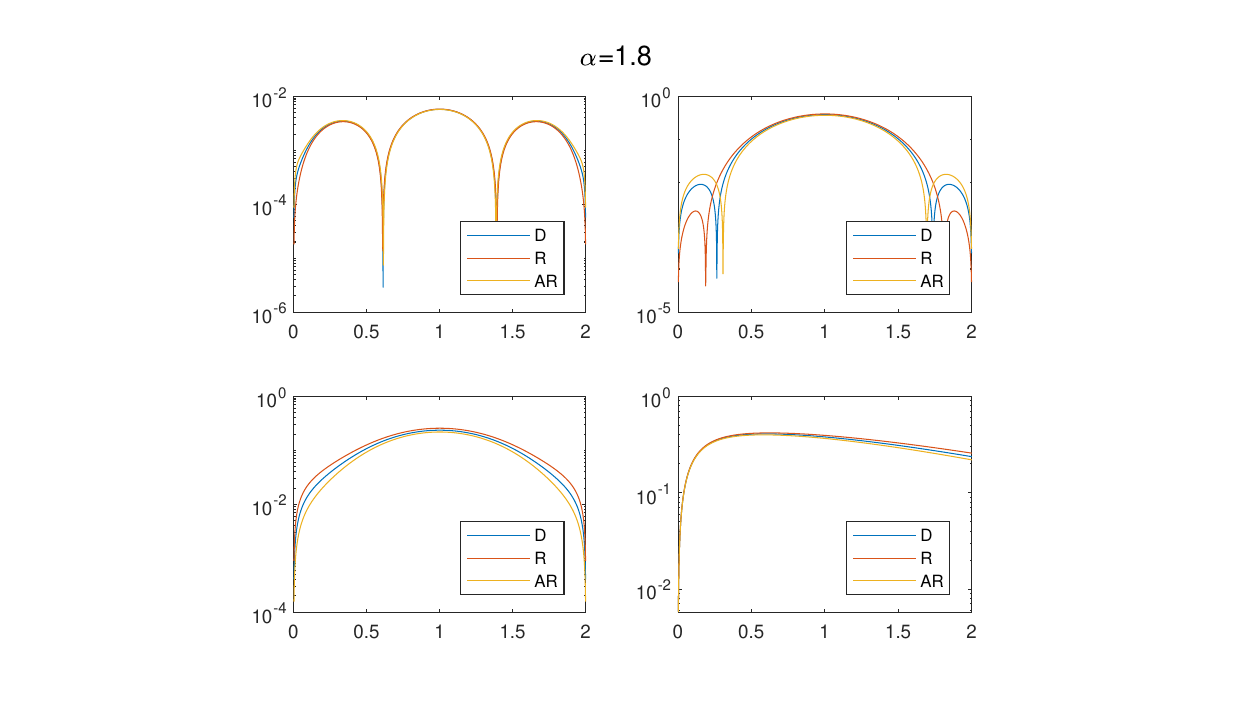}
}
\vskip -1cm \caption{Absolute error vs $x$ at $t \in \{2.e-3,1,2\}$ and Maximal absolute error vs $t$ - Crank-Nicolson,  $\alpha=1.8$ - Homogeneus Dirichlet BCs.}
\label{fig:CN_18x5_omo}
\end{figure}

We consider a test problem  for $\beta=0$ and homogeneous physical BCs
\begin{align}
& \frac{\partial u}{\partial t}(x,t) \label{eq:test}
   =   \Delta_\beta^{\alpha/2} u(x,t)+S(x,t), \quad x\in (0,2), \ t>0, \\ \nonumber
   & u(x,0)= x^4(2-x)^4,  \quad x\in (0,2),\\ \nonumber
   & u(0,t)= 0, \quad t \geq 0, \\ \nonumber
   & u(2,t)= 0, \quad t \geq 0,
\end{align}
with
\[
S(x,t) = e^{-t} \left( -x^4(2-x)^4 -
\frac{1}{2} \sum_{p=0}^4 (-1)^p 2^{4-p} \binom{4}{p} \frac{\Gamma(p+5)}{\Gamma(p+5-\alpha)} (x^{p+4-\alpha} + (2-x)^{p+4-\alpha})\right),
\]
 and exact solution
\[
u(x,t) = e^{-t}x^4(2-x)^4.
\]
As it can be observed in Table \ref{tab:Pb1_bcDo} the numerical anti-reflective BCs guarantees a substantially higher precision for small values of $t$, also in the truncated version. For larger $t$ the difference and advantage are less evident when compared with numerical Dirichlet and reflective BCs. Furthermore, also the parameter $\alpha$ plays a role and for the quality of the reconstruction is better for $\alpha=1.2, 1.4$ while for $\alpha=1.6, 1.8$ the differences are negligible at least for $t=2$. Furthermore, it is worth noticing that in the homogeneous setting numerical anti-symmetric and anti-reflective BCs coincide. In preliminary experiments not reported here, the numerical anti-reflective BCs are definitely better in the non-homogeneous setting: as emphasized in the conclusions, this will be the subject of future research.
Finally, Figures \ref{fig:EI_14x5_omo}, \ref{fig:EI_16x5_omo}, \ref{fig:EI_18x5_omo}, \ref{fig:CN_12x5_omo}, \ref{fig:CN_14x5_omo}, \ref{fig:CN_16x5_omo}, \ref{fig:CN_18x5_omo} reinforce the same conclusions on the higher precision of the numerical anti-symmetric/anti-reflective BCs.

\begin{table}
  \centering
 \footnotesize
  \begin{tabular}{|r|ccc|ccc|}
    \hline
        & \multicolumn {3}{c|}{Implicit Euler } & \multicolumn {3}{c|}{Crank-Nicolson } \\
        \hline
        \multicolumn {7}{|c|}{$\alpha=1.2$}  \\
    \hline
    $t$ & D & R & AR & D & R & AR \\
    \hline
        2.e-03 & 8.931701e-04 & 9.327968e-04 & 8.535608e-04 &  8.927206e-04 & 9.323795e-04 & 8.530704e-04   \\
        1 & 1.932133e-01 & 2.068002e-01 & 1.808326e-01  & 1.929775e-01 & 2.065692e-01 & 1.805912e-01   \\
        2 & 1.782108e-01 & 1.980133e-01  & 1.616497e-01  & 1.778668e-01 & 1.976635e-01 & 1.613089e-01   \\
 %----------------------------------------------------------------------------------------------
    \hline
        \multicolumn {7}{|c|}{$\alpha=1.4$} \\
    \hline
     $t$ & D & R & AR & D & R &  AR \\
    \hline
        2.e-03 & 2.160383e-03 & 2.224234e-03 & 2.096585e-03   & 2.164720e-03 & 2.228610e-03 & 2.100857e-03  \\
        1 & 2.979868e-01 & 3.205744e-01 & 2.779778e-01        & 2.979580e-01 & 3.205535e-01 & 2.779382e-01 \\
        2 & 2.173678e-01 & 2.502798e-01 & 1.916911e-01        & 2.171809e-01 & 2.500784e-01 & 1.915073e-01 \\
    \hline
  %----------------------------------------------------------------------------------------------
    \hline
        \multicolumn {7}{|c|}{$\alpha=1.6$} \\
    \hline
      $t$ & D & R & AR & D & R & AR \\
   \hline
        2.e-03 & 3.804037e-03 & 3.869011e-03 & 3.739125e-03  & 3.819927e-03 & 3.884933e-03 & 3.754953e-03 \\
        1 & 3.474605e-01 & 3.698595e-01 & 3.271661e-01  & 3.474769e-01 & 3.698756e-01 & 3.271757e-01 \\
        2 & 2.326877e-01 & 2.667649e-01 & 2.063672e-01  & 2.326152e-01 & 2.666935e-01 & 2.062872e-01 \\
  %----------------------------------------------------------------------------------------------
    \hline
        \multicolumn {7}{|c|}{$\alpha=1.8$} \\
    \hline
     $t$ & D & R &AR & D & R & AR \\
    \hline
2.e-03 & 5.758829e-03 & 5.800495e-03 & 5.717193e-03  & 5.797546e-03 & 5.839235e-03 & 5.755873e-03  \\
1 & 3.802505e-01 & 3.944576e-01 & 3.668542e-01        & 3.803475e-01 & 3.945544e-01 & 3.669475e-01 \\
2 & 2.372843e-01 & 2.581221e-01 & 2.198008e-01  & 2.372698e-01 & 2.581213e-01 & 2.197737e-01  \\
%----------------------------------------------------------------------------------------------
\hline
  \end{tabular}
  \caption{Maximal absolute error at $t$ in the case of numerical Dirichlet, reflective, and anti-reflective approximations - Homogeneus Dirichlet BCs.} \label{tab:Pb1_bcDo}
\end{table}
%------------------------------------------------

%=================================================================================================================
% Omogenea
%=================================================================================================================
\begin{figure}[h]
{\center
\includegraphics[width=\textwidth]{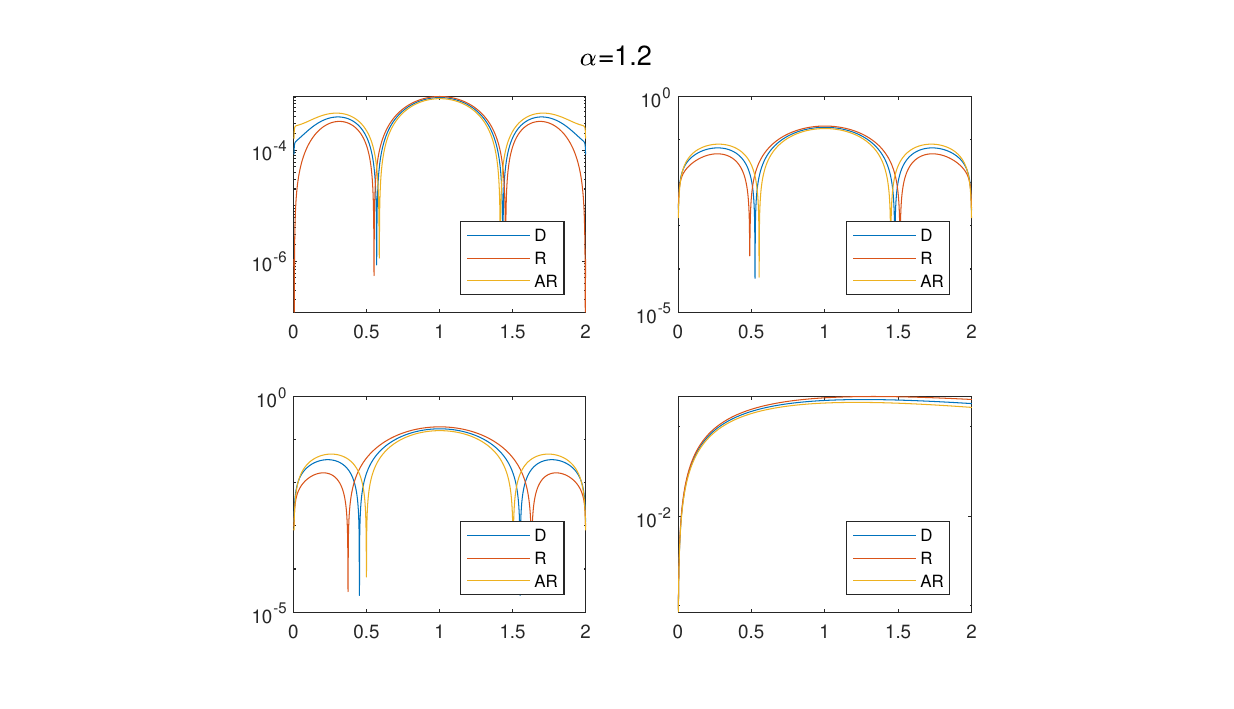}
}
\vskip -1.2cm \caption{Absolute error vs $x$ at $t \in \{2.e-3,1,2\}$ and Maximal absolute error vs $t$ - Implicit Euler,  $\alpha=1.2$ - Homogeneus Dirichlet BCs.}
\label{fig:EI_12x5_omo}
\end{figure}
%--------------------------------------------------CN-------------------------------------------------------------

% ----------------------------------------------------------------
% ----------------------------------------------------------------
\section{Conclusions}\label{sec:end}

In the present work we have combined the idea of fractional differential equations and numerical anti-symmetric and anti-reflective BCs, where the latter  were introduced in a context of signal processing and imaging for increasing the quality of the reconstruction of a blurred signal/image contaminated by noise and for reducing the overall complexity to that of few fast sine transforms i.e. to $O(N\log N)$ real arithmetic operations, where $N$ is the number of pixels. The idea of ending up with matrix structures belonging to the anti-reflective algebra or well approximated by sine transform matrices seems very good also in the considered setting of nonlocal fractional problems.
In fact, we should emphasize that from a matrix viewpoint this is not surprising since both operators, the fractional one and those related to the blurring in imaging are all of nonlocal type. The only relevant difference relies on the subspace of ill-conditioning which corresponds to low frequencies in the fractional differential case and to the high frequencies in the case of blurring.

Several numerical tests, tables, and visualizations have been provided and critically discussed, also in connection with other classical numerical BCs and the truncation of the two types of numerical BCs mainly considered in the current work (anti-reflective and anti-symmetric).

More research remains to be performed especially in connection with the following items:
\begin{itemize}
\item multidimensional domains and more involved nonlocal operators which could be treated since the $\tau$ algebra and the related sine I and anti-reflective transforms admit multilevel versions via tensorizations \cite{S-AR-proposal,paperART};
\item of course, if the considered matrices do not belong to an algebra and preconditioned Krylov solvers have to be employed, then it is necessary to take into account the computational barriers in \cite{nega1,nega2};
\item non-equispaced grids or variable coefficients which could be treated thanks to the theory of generalized locally Toeplitz matrix-sequences \cite{GLT-exposition-eng} also in a multidimensional setting \cite{GLT-BookII,Doro};
\item a numerical and theoretical comparison in terms of precision between the truncated and nontruncated versions of the considered numerical BCs, even if few numerics suggest that the differences are not relevant;
\item a delicate issue concerns the case of non-homogeneous physical Dirichlet BCs. Non-homogeneous BCs lead to a further term in the operator and to a corresponding further term in its numerical approximation \cite{treat-NonHom1,treat-NonHom2}, so that the relations (\ref{glapproxl}), (\ref{glapproxr}) are no longer shift invariant. In other words, the rectangular Toeplitz structures are lost and the associated matrices have an additional term, which is still Toeplitz-like but not purely Toeplitz. As a consequence the way of adding the most precise numerical BCs is nontrivial and it will be the subject of future investigations.
\end{itemize}

The last two items are of particular interest since they have the potential to give a formal answer to the problem of unphysical oscillations (rough artifacts similar to those in imaging \cite{imaging-reflective,S-AR-proposal}) appearing in the numerical solution, close to the boundaries and for small stepsizes, as observed in \cite{BU2014,rough-artifact} for stationary fractional problems.

\section*{Acknowledgements}

Ercília Sousa was partially supported  by the Centre for Mathematics of the University of Coimbra – UIDB/00324/ 2020, funded by the Portuguese Government through FCT/MCTES.
The work of Stefano Serra-Capizzano and Cristina Tablino-Possio is supported by GNCS-INdAM.
The work of Rolf Krause and Stefano Serra-Capizzano is funded from the European High-Performance Computing Joint Undertaking  (JU) under grant agreement No 955701. The JU receives support from the European Union’s Horizon 2020 research and innovation programme and Belgium, France, Germany, Switzerland.
Furthermore, Stefano Serra-Capizzano is grateful for the support of the Laboratory of Theory, Economics and Systems – Department of Computer Science at Athens University of Economics and Business.

\end{document}